\documentclass[11pt]{article}
\usepackage{amscd,amssymb,stmaryrd}
\usepackage[fleqn]{amsmath}
\usepackage{graphicx}
\usepackage{subcaption}
\usepackage[utf8]{inputenc}
\usepackage[export]{adjustbox}
\usepackage{wrapfig}
\usepackage{amstext}
\usepackage{multirow}
\usepackage{pstricks,pst-node,pst-plot,pst-coil}
\usepackage{amsthm}
\usepackage{algpseudocode}
\usepackage{algorithm}
\usepackage[title]{appendix}
\usepackage{amsmath}
\usepackage{color}
\usepackage{mathtools}
\usepackage{comment}
\usepackage{hyperref}
\usepackage[english]{babel}
\usepackage{booktabs}
\usepackage{mathrsfs}
\usepackage{tcolorbox}

\newcommand\norm[1]{\left\|#1\right\|}
\newcommand\abs[1]{\lvert#1\rvert}

\newcommand{\tand}{\quad \text{and} \quad}
\newcommand{\mR}{\mathcal{R}}

\newcommand{\eps}{\varepsilon}
\newcommand{\revi}[1]{\textcolor{black}{#1}}
\newcommand{\revii}[1]{\textcolor{black}{#1}}

\theoremstyle{definition}
\newtheorem{theorem}{Theorem}[section]

\newtheorem{example}[theorem]{Example}
\theoremstyle{remark}

\title{\sc Temporal Splitting algorithms for non-stationary multiscale problems}
\author{Yalchin Efendiev\footnote{Department of Mathematics, Texas A\&M University, College Station, TX 77843, USA \& North-Eastern Federal University, Yakutsk, Russia}, ~ Sai-Mang Pun\footnote{Department of Mathematics, Texas A\&M University, College Station, TX 77843, USA}, ~ Petr N. Vabishchevich\footnote{Nuclear Safety Institute, Russian Academy of Sciences, Moscow, Russia \& North-Eastern Federal University, Yakutsk, Russia}}
\begin{document}
\maketitle
\begin{abstract}
In this paper, we study temporal splitting algorithms for multiscale 
problems. 
The exact fine-grid spatial problems typically require some
reduction in degrees of freedom. Multiscale algorithms
are designed to represent the fine-scale details on a coarse grid
and, thus, reduce the problems' size. When solving time-dependent problems,
one can take advantage of the multiscale decomposition of the solution
and perform temporal splitting by solving smaller-dimensional problems,
which is studied in the paper. In the proposed approach, we consider
the temporal splitting based on various low dimensional spatial approximations.
Because a multiscale spatial splitting gives a ``good'' decomposition of
the solution space, one can achieve an efficient implicit-explicit temporal discretization. 
We present a recently developed theoretical result in \cite{efendiev2020splitting} and adopt in this paper for multiscale problems. Numerical results are presented to demonstrate the efficiency of the proposed splitting algorithm. 

\end{abstract}


\section{Introduction}

Multiscale phenomena occur in many applications such as 
porous media, material sciences, and biological applications. 
Detailed simulations of multiscale problems are prohibitively
expensive due to the fact that one needs to resolve the smallest scales. 
For this reason, various multiscale algorithms have been developed
to reduce the problem's size. Typically, these algorithms
are used to reduce the spatial dimensions. Some space-time
multiscale methods are developed to deal with spatial and temporal scales.
In many cases, appropriate temporal splitting can be used to alleviate
the computational complexity and this is the goal of this paper.

Multiscale algorithms have been intensively developed for handling
spatial multiscale features \cite{chung2016adaptiveJCP,eh09}. 
These algorithms typically use multiscale
basis functions or some derived effective properties to reduce
the spatial complexity of the problem. Meanwhile, splitting
algorithms \cite{glowinski2017splitting,marchuk1990splitting,VabishchevichAdditive} have been extensively studied for non-stationary problems,
which are decomposed into several smaller sub-problems and
solved in sub-time steps. Combining spatial reduction and temporal splitting can reduce the computational cost associated to the large-scale problem.

Multiscale spatial
algorithms have been extensively studied in the literature.
These include
homogenization-based approaches \cite{eh09,le2014msfem}, 
multiscale
finite element methods (MsFEM) \cite{eh09,hw97,jennylt03}, 
generalized multiscale finite element methods (GMsFEM) \cite{chung2016adaptiveJCP,MixedGMsFEM,WaveGMsFEM,chung2018fast,GMsFEM13}, constraint energy minimizing GMsFEM (CEM-GMsFEM) \cite{chung2018constraint, chung2018constraintmixed}, nonlocal
multi-continua (NLMC) approaches \cite{NLMC},
metric-based upscaling \cite{oz06_1}, heterogeneous multiscale method \cite{ee03}, localized orthogonal decomposition (LOD) \cite{henning2012localized}, equation-free approaches \cite{rk07,skr06}, 
and hierarchical multiscale method \cite{brown2013efficient}. 
Some of these approaches, such as homogenization-based
approaches, are designed for problems with scale 
separation, which is a reasonable assumption in the modeling of novel material. In porous media applications, the spatial heterogeneities
are too complex to have scale separation. In addition, they 
contain large jumps in the coefficients of the model. 
For these purposes, 
a general concept is introduced for CEM-GMsFEM and NLMC approaches, where multiple basis
functions or continua are designed to capture the multiscale features well and solve problems on a coarse grid
\cite{chung2018constraintmixed, NLMC}, which does not necessarily resolve all possible scales and contrast. 
These approaches require a careful design of multiscale
dominant modes. We remark that the application of these methods to hyperbolic
equations i challenging \cite{WaveGMsFEM,chung2020computational} due to distant temporal effects.


Standard splitting schemes are based on the initial additive 
decomposition of the differential operator
\cite{marchuk1990splitting,VabishchevichAdditive}. 
Typical examples include convection-diffusion, where convection and diffusion are used in the splitting. In these cases, the operator is decomposed
based on physical processes.
In many cases, it is more appropriate to consider solution decomposition, 
which indirectly yields a decomposition of the operator. 
This is the case for multiscale decomposition, where 
the solution space is divided into several coarse-grid components.
In this  work,
we will make use of a recent decomposition developed in 
 \cite{efendiev2020splitting}, where a general framework is proposed and 
the transition to simpler problems is carried out 
based on spatial decomposition of the operator.
Based on this decomposition, we construct splitting schemes for the solution. 
An approximate solution to the Cauchy problem in a finite-dimensional 
Hilbert space for a first-order evolution equation is constructed 
on a family of spaces using proposed restriction and prolongation operators.
After that, the individual components of the solution are determined from the system of evolutionary equations.
Three-level schemes for splitting the solution with time transfer to the upper level are proposed and investigated based on splitting that uses
the diagonal parts of the operator matrices of the corresponding evolutionary system for the components of the solution. 

In the current paper,
this general framework is used to solve decompose space-time features
of the multiscale solution. 
We divide the space into various multiscale
subspaces and couple them via splitting algorithms.
In particular, we combine the temporal splitting algorithms and spatial
multiscale methods. We divide the spatial space into various components
and use these subspaces in the temporal splitting. As a result, smaller
systems are inverted in each time step, which reduces the computational
cost. Our main goal is (i) to reduce the computational cost of multiscale spatial algorithms for time-dependent problems and 
(ii) to set up a framework for studying the interplay between
 temporal and spatial scales. In implicit methods,
the resulting algorithms give smaller dimensional problems to invert.
 We formulate our main theoretical result and show 
that the decomposition method is consistent and unconditionally stable provided that appropriate assumption for user-defined parameters is satisfied. 
Numerical results show that one can achieve high accuracy at a lower
computational cost. In our numerical results, we compare the results
of the proposed algorithms with backward Euler scheme. 

The paper is organized as follows. In Section \ref{sec:model-setting}, we present some preliminaries and the model problem considered in this work. 
Section \ref{sec:gmsfem} is devoted to the multiscale methods and we will briefly overview the construction of multiscale basis function within the framework of the GMsFEM. In Section \ref{sec:splitting}, we present
our solution decomposition algorithm and provide a theoretical result. In Section \ref{sec:numerics}, we present
some numerical results to demonstrate the efficiency of the proposed method. Concluding remarks are drawn in Section \ref{sec:conclusion}. 

\section{Preliminaries} \label{sec:model-setting}

In this section, we present the model problem considered in this work. 
We consider the following parabolic partial differential equation (PDE) defined in the computational domain $\revi{[0,T]\times \Omega}$ (with $T>0$) as follows: find a function $U$ such that 
\begin{eqnarray}
\begin{split}
\frac{\partial U}{\partial t} - \nabla \cdot (\kappa \nabla U) &= \mathcal{F} &\quad \text{for} ~ \revi{(t,x) \in (0,T] \times \Omega},\\
U\revi{(t,x)} & = 0 &\quad \text{for} ~ \revi{(t,x) \in (0,T] \times \partial \Omega},\\
U\revi{(0,x)} & = U_0(x) &\quad \text{for} ~ x \in \Omega. 
\end{split}
\label{eqn:model-pde}
\end{eqnarray}
Here, the function $\kappa \in L^\infty (\Omega)$ is a heterogeneous permeability field. We assume that there exist two constants $0< \kappa_{\min} \ll \kappa_{\max}$ such that 
$$\kappa_{\min} \leq \kappa(x) \leq \kappa_{\max} \quad \text{for almost every} ~ x \in \Omega.$$
The function $\mathcal{F} \in L^2([0,T]; L^2(\Omega))$ is a source term and $U_0 \in L^2(\Omega)$ is a given initial condition. For simplicity, we assume that the solution $U$ satisfy the homogeneous boundary condition $U = 0$ on the boundary $\partial \Omega$. The extension to non-homogeneous boundary condition is straightforward. 

In this work, we consider the discretization of the problem \eqref{eqn:model-pde} using the generalized multiscale finite element method (GMsFEM). The framework of GMsFEM provides a systematical approach to generate multiple degrees of freedom per node from a coarse grid of the computational domain. This set of coarse degrees of freedom yields a natural option for the solution decomposition algorithm derived in the previous work \cite{efendiev2020splitting}. In the following sections, we will introduce the framework of GMsFEM and derive the solution decomposition scheme. 
\revi{Using the GMsFEM for the spatial discretization}, we then arrive at the following Cauchy problem for the first-order evolution: find $u = u(t) \in \mathbb{R}^{\mathcal{N}_c}$ such that 
\begin{eqnarray}
\begin{split}
M \frac{du}{dt} + Au & = f  \quad \text{for} ~ t \in (0,T], \\
u(0) & = u_0,
\end{split}
\label{eqn:model}
\end{eqnarray}
with $M$ and $A$ being the mass and stiffness matrices, respectively. 
That is, the matrices $M$ and $A$ are defined by 
$$ M := \left ( (\varphi_j, \varphi_i)_\Omega \right ) \in \mathbb{R}^{\mathcal{N}_{\text{c}} \times \mathcal{N}_{\text{c}}} \tand A := \left ( a(\varphi_j, \varphi_i) \right ) \in \mathbb{R}^{\mathcal{N}_{\text{c}} \times \mathcal{N}_{\text{c}}},$$
where \revi{$(v, w)_\Omega := \int_\Omega vw ~ dx$ is the standard $L^2$ inner product in the domain $\Omega$ for any function $v, w \in L^2(\Omega)$, $a(v,w) := \int_\Omega \kappa \nabla v \cdot \nabla w ~ dx$ for any $v, w \in H_0^1(\Omega)$}, and $\{ \varphi_i \}_{i=1}^{\mathcal{N}_{\text{c}}}$ is a set of multiscale basis functions with cardinality $\mathcal{N}_{\text{c}} \in \mathbb{N}^+$. 
We will briefly review the construction of the multiscale basis functions in Section \ref{sec:gmsfem}. 
The right-hand side $f \in \mathbb{R}^{\mathcal{N}_{\text{c}}}$ and the initial condition $u_0 \in \mathbb{R}^{\mathcal{N}_{\text{c}}}$ in \eqref{eqn:model} are simply defined by 
$$ f := \left ( \begin{array}{c}
(\mathcal{F}, \varphi_1)_\Omega \\
\vdots \\
(\mathcal{F}, \varphi_{\mathcal{N}_{\text{c}}})_\Omega \\
\end{array} \right ) \tand u_0 := \left ( \begin{array}{c}
(U_0, \varphi_1)_\Omega \\
\vdots \\
(U_0, \varphi_{\mathcal{N}_{\text{c}}})_\Omega \\
\end{array} \right ).$$
\revi{In the GMsFEM setting, the solution vector $u \in \mathbb{R}^{\mathcal{N}_c}$ represents the vector of coefficients with respect to the set of multiscale basis functions. }
We remark that the problem \eqref{eqn:model} also arises in other spatial discretization of PDEs. For instance, when using finite difference scheme to discretize parabolic PDEs, \revi{the solution of \eqref{eqn:model}} is a grid-based function defined at the nodes of a computational grid of the domain. 
In general, one considers the equation \eqref{eqn:model} over a finite-dimensional Hilbert space $\mathcal{U}$ with prescribed right-hand side $f \in \mathcal{U}$ and the initial condition $u_0 \in \mathcal{U}$. 
\revi{In this case, the operators $M: \mathcal{U} \to \mathcal{U}$ and $A: \mathcal{U} \to \mathcal{U}$ are self-adjoint and positive definite. One may assume that the operators $M$ and $A$ do not depend on the temporal variable.} 

\section{The multiscale method} \label{sec:gmsfem}
In this section, we briefly overview the framework of GMsFEM. 
For further details on GMsFEM we refer the reader to
\cite{Review12,GMsFEM13,eglp13}, 
and the references therein.
The framework of this systematic approach starts with the construction of snapshot functions. 
After that, one may obtain the multiscale basis functions by solving a class of specific spectral problems in the snapshot space, and these multiscale basis functions will be used to solve the multiscale solution. 
As an example, we consider the following elliptic boundary-value problem in a spatial domain $\Omega \subset \mathbb{R}^d$ 
\begin{eqnarray}
\mathcal{L} (U)  = g \quad \text{in } \Omega \tand U = 0 \quad \text{on } \partial \Omega
\label{eqn:elliptic}
\end{eqnarray}
in order to illustrate the main framework of the GMsFEM, where the elliptic differential operator $\mathcal{L}(\cdot) := - \nabla \cdot (\kappa \nabla (\cdot) )$ contains some multiscale features. 

We start with the notion of fine and coarse grids. Let $\mathcal{T}^H$ 
be a conforming partition of the computational domain $\Omega$ with mesh size $H>0$.
We refer to this partition as the coarse grid.
Subordinating to the coarse grid, we define a fine grid 
partition (with mesh size $h \ll H$), denoted by $\mathcal{T}^h$, by refining each coarse element $K \in \mathcal{T}^H$ into a connected union of fine-grid blocks. We assume the above refinement is performed such that $\mathcal{T}^h$ is a conforming partition of the domain $\Omega$. 
Let $\{x_i\}_{i=1}^{N_c}$ 
be the set of coarse-grid nodes of the coarse grid $\mathcal{T}^H$ with cardinality $N_c$. 
Let $N$ be the number of elements in the coarse grid. 
We define coarse neighborhood of the node $x_i$ by 
$$\omega_i := \bigcup \{ K_j \in \mathcal{T}^H: x_i \in \overline{K}_j \};$$
that is, the union of all coarse elements which have the node $x_i$ as a vertex. We call $\omega_i$ the coarse neighborhood corresponding to the node $x_i$.

\subsection{Snapshot space}
We first present the construction of the snapshot space. The snapshot space is pre-constructed before solving the actual problem. The snapshot space consists of harmonic extensions of fine-grid functions that are defined on the generic neighborhoods $\omega_i$. 

We define the fine-grid function $\delta_\ell^h(x_k) := \delta_{\ell k}$
for $x_k \in J_h(\omega_i)$, where $J_h(\omega_i)$ denotes the set of fine-grid boundary nodes on $\partial \omega_i$. 
Here, $\delta_{\ell k}$ is the Kronecker symbol.
Denote $L_i$ the cardinality of $J_h(\omega_i)$. Then, for $\ell= 1,\cdots,L_i$, the snapshot function $\eta_\ell^{(i)}$ is defined to be the solution to the following system
\begin{eqnarray*}
    \mathcal{L}^\eps (\eta_\ell^{(i)})  = &  0 &\quad \text{in } K \subset \omega_i, \\
    \eta_{\ell}^{(i)} = & g_i & \quad \text{on } \partial K \text{ for each } K \subset \omega_i,\\
    \eta_{\ell}^{(i)}  = & \delta_\ell^h & \quad \text{on } \partial \omega_i.
\end{eqnarray*}
The function $g_i$ is linear and continuous along the boundary $\partial K$ for each coarse element $K \in \mathcal{T}^H$ satisfying $K \subset \omega_i$. The local snapshot space $V_{\text{snap}}^{(i)}$ corresponding to the coarse neighborhood $\omega_i$ is defined as $ V_{\text{snap}}^{(i)} := \text{span} \{ \eta_\ell^{(i)}\}_{\ell = 1}^{L_i}$. One may define the global snapshot space $V_{\text{snap}}$ as $ V_{\text{snap}} := \bigoplus_{i=1}^{N_c} V_{\text{snap}}^{(i)}.$

\subsection{Offline multiscale basis construction}
In order to reduce the number of degrees of freedom in the snapshot space, we perform a spectral decomposition in the snapshot space and select the dominant modes (corresponding to small eigenvalues) to construct the multiscale space. 
As a result, we obtain a multiscale space that will be used in the simulation. 

For each $i = 1,\cdots, N_c$, the spectral problem reads as follows: find $\phi_j^{(i)} \in V_{\text{snap}}^{(i)}$ and $\lambda_j^{(i)} \in \mathbb{R}$ such that
\begin{eqnarray}
 a_i(\phi_j^{(i)}, w) = \lambda_j^{(i)} s_i(\phi_j^{(i)}, w) \quad \forall w \in V_{\text{snap}}^{(i)}, \quad j = 1,\cdots, L_i,
 \label{eqn:sp}
\end{eqnarray}
where $a_i(\cdot,\cdot)$ is a symmetric non-negative definite bilinear form and $s_i(\cdot,\cdot)$ is a 
symmetric positive definite bilinear operators defined on 
$V_{\text{snap}}^{(i)} \times V_{\text{snap}}^{(i)}$. The eigenfunctions  
$\phi_j^{(i)}$ are normalized to satisfy $s_i(\phi_j^{(i)},\phi_j^{(i)})=1$. 
Assume that the eigenvalues obtained from \eqref{eqn:sp} are arranged in ascending order 
and we use the first $\ell_i \in \mathbb{N}^+$ eigenfunctions (corresponding to the smallest $\ell_i \leq L_i$ 
eigenvalues) to construct local auxiliary multiscale space 
$$V_{\text{off}}^{(i)} := \text{span} \left \{\psi_j^{(i)} := \chi_i \phi_j^{(i)} | \  j= 1,\cdots, \ell_i \right \}.$$
Here, we multiply the eigenfunctions $\phi_j^{(i)}$ by the partition of unity function $\chi_i$ related to the coarse neighborhood $\omega_i$ to obtain a set of conforming basis functions. The partition of unity function $\chi_i$ is defined on each coarse neighborhood $\omega_i$ such that $\chi_i = 0$ on $\partial \omega_i$ and $\sum_{i=1}^{N_c} \chi_i = 1$ (see \cite{hw97} for more details about the standard multiscale partition of unity function). 
The global auxiliary space $V_{\text{off}}$ is the direct sum of these local auxiliary 
multiscale space, namely $V_{\text{off}} := \bigoplus_{i=1}^{N_c} V_{\text{off}}^{(i)}$ with dimension $\text{dim}(V_{\text{off}}) =: \mathcal{N}_{\text{c}} = \sum_{i=1}^{N_c} \ell_i$. 

Notice that each basis function in $V_{\text{off}}$ is represented on the fine grid $\mathcal{T}^h$ by a set of  basis functions in a fine-scale finite element space $V_h$ with dimension $\text{dim}(V_h) =: \mathcal{N}_{\text{f}} \gg \mathcal{N}_{\text{c}}$. Therefore, each $\psi_j^{(i)}$ can be represented by a vector $\Psi_j^{(i)}$ containing the coefficients in the expansion in the fine-grid basis functions. We can define prolongation matrix $\mathcal{R}_{\text{off}}$ such that 
$$ \mathcal{R}_{\text{off}} := \left ( \Psi_1^{(1)}\  \cdots \ \Psi_{\ell_{N_c}}^{(N_c)} \right ) \in \mathbb{R}^{\mathcal{N}_{\text{f}} \times \mathcal{N}_{\text{c}}},$$
which maps from the coarse component to the fine one. 
The construction of multiscale basis functions provides a natural option to decompose the multiscale solution. We will perform solution decomposition based on the structure of the prolongation matrix. 

We remark that the choice of $a_i(\cdot,\cdot)$ and $s_i(\cdot,\cdot)$ is problem-dependent and motivated by the analysis of the problem considered. In general, the bilinear forms $a_i(\cdot,\cdot)$ and $s_i(\cdot,\cdot)$ are chosen such that for any $\psi \in V_{\text{snap}}^{(i)}$ there exist a function $\psi_0 \in V_{\text{off}}^{(i)}$ and a generic constant $\delta >0$, which is independent of the heterogeneities, satisfying
$$ a_i(\psi - \psi_0, \psi - \psi_0) \leq \delta s_i(\psi - \psi_0,\psi - \psi_0).$$

\section{The solution decomposition algorithm} \label{sec:splitting}
In this section, we derive the solution decomposition algorithm. 
Following \cite{efendiev2020splitting}, we assume that the solution to \eqref{eqn:model} has an additive representation as follows: 
\begin{eqnarray}
u (t) = \sum_{i=1}^p u_i(t)
\label{eqn:sol-split}
\end{eqnarray}
for some positive integer $p \geq 2$. 
Then, we construct each individual term $u_i$ to seek the solution and achieve the goal of simplifying the tasks for finding these terms. 
To simplify the discussion below, we assume that $p = 2$ and divide the solution into two parts. The extension of the case with larger $p$ is straightforward. 

In this work, we assume that each component in the additive representation is spanned by different basis functions in the multiscale space $V_{\text{off}}$. The construction of multiscale basis functions provides a natural option to decompose the multiscale solution according to that of the basis functions. As a result, 
we compute the first part of the solution $u_1$ using first set of dominant basis functions while the second part $u_2$ is computed using another set of dominant multiscale basis functions. We consider the decomposition of  (the coefficients of) multiscale solution. Let $u \in \mathbb{R}^{\mathcal{N}_{\text{c}}}$ be the solution of \eqref{eqn:model} and it represents the vector of coefficients with respect to the set of multiscale basis functions. 
Specifically, one may write 
\begin{eqnarray}
 \mathcal{R}_{\text{off}} (u) = \mathcal{R}_1 V + \mathcal{R}_2 W,
\label{eqn:sol-split-1}
\end{eqnarray}
where $\mathcal{R}_1$ and $\mathcal{R}_2$ are $\mathcal{N}_{\text{f}} \times m$ and $\mathcal{N}_{\text{f}} \times k$ (prolongation) matrices with $m+k = \mathcal{N}_{\text{c}}$, $V \in \mathbb{R}^m$, and $W \in \mathbb{R}^k$. For instance, up to a permutation of the columns, one may set $\mathcal{R}_1$ and $\mathcal{R}_2$ such that 
$$ \mathcal{R}_{\text{off}} = \left ( \mathcal{R}_1 ~ \mathcal{R}_2 \right ).$$
Substituting the formula \eqref{eqn:sol-split-1} to \eqref{eqn:model} and multiplying by the {\it test spaces} $\mR_1^T$ and $\mR_2^T$, we obtain 
\begin{eqnarray}
\begin{split}
\mR_1^T M \mR_1 \frac{dV}{dt} + \mR_1^TM \mR_2 \frac{dW}{dt} + \mR_1^T A \mR_1 V + \mR_1^T A \mR_2 W & = \mR_1^T f, \\
\mR_2^T M\mR_1 \frac{dV}{dt} + \mR_2^TM \mR_2 \frac{dW}{dt} + \mR_2^T A \mR_1 V + \mR_2^T A \mR_2 W & = \mR_2^T f, \\
\end{split}
\label{eqn:model-1}
\end{eqnarray}
or equivalently 
\begin{eqnarray}
C \frac{d\revi{Z}}{dt} + B\revi{Z} = \tilde f,
\label{eqn:model-2}
\end{eqnarray}
where 
$$ C = \left ( \begin{array}{cc}
C_{11} & C_{12} \\
C_{21} & C_{22}
\end{array} \right ) = \left ( \begin{array}{cc}
\mR_1^TM \mR_1 & \mR_1^TM \mR_2 \\
\mR_2^TM \mR_1 & \mR_2^TM \mR_2 
\end{array} \right ), \quad B = \left ( \begin{array}{cc}
B_{11} & B_{12} \\
B_{21} & B_{22}
\end{array} \right ) = \left ( \begin{array}{cc}
\mR_1^T A \mR_1 & \mR_1^T A \mR_2 \\
\mR_2^T A \mR_1 & \mR_2^T A \mR_2 
\end{array} \right ), $$
$$
\revi{Z} = \left ( \begin{array}{c}
V \\ W \end{array} \right ), \tand \tilde f = \left ( \begin{array}{c}
\mR_1^T f \\
\mR_2^T f \end{array} \right ).
$$
We remark that in the setting of GMsFEM, the matrices $C$ and $B$ are 
the mass and stiffness matrices defined on the coarse grid. Therefore, these matrices are symmetric and positive definite. 
We will consider the splitting of the system based on multiscale decomposition
of the space. With the proposed splitting,  we reduce
 the dimension of the
system that arise in implicit discretization, and also associate separate space
and time scales.

Next, we perform the temporal discretization for the system \eqref{eqn:model-2}. 
We choose a time step $\tau>0$ such that $\tau = T/N$ with $N \in \mathbb{N}^+$.  
We assume that the matrices $C$ and $B$ can be written as $C = C_1 + C_2$ and $B = B_1 + B_2$ and 
consider the following temporal splitting scheme of \eqref{eqn:model-2}: given $Z^0$ and $Z^1$, find $\{ Z^{n+1} \}_{n=1}^{N-1}$ such that 
\begin{eqnarray}
\begin{split}
C_1 \left ( \mu \frac{Z^{n+1} - Z^n}{\tau} + (1-\mu) \frac{Z^{n} - Z^{n-1}}{\tau} \right ) &+ C_2 \frac{Z^{n} - Z^{n-1}}{\tau} \\
& + B_1 (\sigma Z^{n+1} + (1-\sigma)Z^n) + B_2 Z^n = \tilde f^{n+1}, 
\end{split}
\label{eqn:splitting-scheme}
\end{eqnarray}
where $\tilde f^n = \tilde f(t_n)$ and $t_n = n\tau$ for $n = 1,\cdots, N-1$. 
The term $Z^0$ is determined by the original initial condition $u_0$ while the term $Z^1$ can be computed, for instance, by the backward Euler scheme as follows: 
\begin{eqnarray}
C \frac{Z^1 - Z^0}{\tau} + BZ^1 = \tilde f^1. 
\label{eqn:init_z1}
\end{eqnarray}
One can select the block-diagonal parts of $C$ and $B$ to form the matrices $C_1$ and $B_1$. That is, we can set 
$$ C_1 = \left ( \begin{array}{cc} C_{11} & O \\
O & C_{22} \end{array} \right ) \quad \text{and} \quad 
B_1 = \left ( \begin{array}{cc} B_{11} & O \\
O & B_{22} \end{array} \right ).$$
Here, we denotes $O$ the zero matrix to match the size of the blocks. Similarly, $C_1$ and $B_1$ are symmetric and positive definite. We remark that there is other possible option to split the matrices. For example, we can choose 
$$ C_1 = \left ( \begin{array}{cc} C_{11}/2 & O \\
C_{21} & C_{22}/2 \end{array} \right ) \quad \text{and} \quad 
B_1 = \left ( \begin{array}{cc} B_{11}/2 & O \\
B_{21} & B_{22}/2 \end{array} \right ).$$
Then, we have $C_2 = C_1^T$ and $B_2 = B_1^T$. We mainly focus on the diagonal splitting in the section  of numerical experiments. 

The parameters $\mu$ and $\sigma$ are introduced to balance the solution decomposition. One may need to choose the parameters appropriately in order to get an unconditionally stable numerical scheme. To find an approximated solution on a new level, we solve the following system:
$$ (\mu C_1 + \tau \sigma B_1) Z^{n+1} = \varphi^n,$$
where 
\begin{eqnarray*}
 \varphi^n = \tau \tilde f^{n+1} - \Big ( \tau (1-\sigma) B_1 + \tau B_2 + (1-2\mu)C_1 + C_2 \Big ) Z^n + \Big (  (1-\mu) C_1 + C_2 \Big ) Z^{n-1}. 
\end{eqnarray*}
Next, we present the stability of the proposed splitting method and this can be proved by adopting the techniques in \cite[Theorem 1]{efendiev2020splitting}.
For any matrix $A$, we write $A\geq 0$ if $A$ is non-negative definite and $A >0$ if $A$ is positive definite. 
Then, one can show that the scheme \eqref{eqn:splitting-scheme}-\eqref{eqn:init_z1} is consistent. Moreover, assume that $B_1$ is non-negative definite; this can be guaranteed if diagonal splitting is performed.  
If the parameters $\mu>0$ and $\sigma>0$ satisfy 
\begin{eqnarray}\label{eqn:sigma-mu}
\mu C_1 - \frac{1}{2} C >0 \tand  \sigma  B_1 - \frac{1}{4} B >0,
\end{eqnarray}
the temporal splitting scheme \eqref{eqn:splitting-scheme}-\eqref{eqn:init_z1} can be shown to be unconditionally stable. Furthermore, the following a priori estimate holds: 
$$ \norm{\frac{Z^{n+1} + Z^{n}}{2}}_B^2 \leq \frac{1}{\tau^2} \norm{Z^1 - Z^0}_D^2 + \norm{\frac{Z^1 + Z^0}{2}}_B^2 + \frac{1}{2} \sum_{i=1}^n \tau \norm{\tilde f^{i+1}}^2,$$
where $\norm{\cdot}$ denotes both the Euclidean norm in $\mathbb{R}^{\mathcal{N}_{\text{c}}}$. 
\revi{The norms $\norm{\cdot}_B$ and $\norm{\cdot}_D$ are defined by 
\begin{eqnarray*}
\begin{split}
\norm{v}_B &:= \sqrt{v^T B v}, \\
\norm{v}_D &:= \sqrt{v^T Dv} \quad \text{with} \quad D := \tau \left ( \mu C_1 - \frac{1}{2} C \right ) + \tau^2 \left ( \sigma B_1 - \frac{1}{4} B \right ),
\end{split}
\end{eqnarray*}
for any vector $v \in \mathbb{R}^{\mathcal{N}_c}$. 
These norms are well-defined since $B$ is symmetric and positive definite; $D$ is also symmetric and positive definite if \eqref{eqn:sigma-mu} is satisfied. 
}
We remark that by Lax-Richtmyer equivalence theorem \cite{lax1956survey} (see also Sections 10.5 and 10.6 in \cite{strikwerda2004finite}), the stability result above implies that the discrete solution $Z^n$ converges to the solution of the continuous problem. 

\revi{The proposed splitting method provides a framework to split the original discretized system \eqref{eqn:model-2} into a block-diagonal three-level system \eqref{eqn:splitting-scheme} to reduce the computational cost. The proposed splitting also allows one to deal with fast and slow temporal scales separately. This can be used to design implicit-explicit schemes, treat nonlinearities, or stiff terms appropriately.} 

\revi{We remark that the selection of splitting and the temporal step size are fixed in the numerical experiments below. We remark that one can adjust splitting during the time marching process as more basis functions will be required to capture the dynamics. Also, one can design adaptive time marching scheme combining the proposed splitting technique to make the algorithm more efficient. We will study these issues in the future.} 

We briefly discuss the error between the splitting solution $\{ Z^n \}$ and the one $\{ \tilde Z^n \}$ obtained by the backward Euler method. Denote $e^n := \tilde Z^n - Z^n$. To simplify the discussion, we assume $\mu = \sigma = 1$. Then, in this case, we have 
$$ (C+ \tau B) \revii{e^{n+1}} = C e^n + C_2 (Z^n - Z^{n-1}) - (C_2 + \tau B_2) (Z^{n+1} - Z^n).$$
In the setting of diagonal splitting, the error between the splitting algorithm and the one obtained by the backward Euler method depends on the structure of the off-diagonal parts $C_2 + \tau B_2$. 
\revii{From the formula above, the quantity $\norm{e^{n+1}}$ is small if the singular values of the non-zero block entries of the matrix $C_2 + \tau B_2$ becomes small.} To achieve this, one can perform orthogonalization such that most of the basis functions in $\mathcal{R}_1$ are orthogonal to that of $\mathcal{R}_2$. 
In our numerical experiments, we orthogonalize multiscale basis functions
in each coarse region. This reduces the effect of the off-diagonal matrices 
in the forward simulations.
Another mean of reducing off-diagonal terms is possibly using Constraint
Energy Minimizing GMsFEM (CEM-GMsFEM) (proposed originally in \cite{chung2018constraint}).
We will consider this issue in our future works.

Furthermore, for any positive integer $p\geq 2$ in the additive representation \eqref{eqn:sol-split}, one can also show that under the condition on the parameters 
$$ \mu \geq \frac{p}{2} \quad \text{and} \quad \sigma \geq \frac{p}{4},$$
the block-diagonal splitting scheme is unconditionally stable with the matrices $C_1$ and $B_1$ being 
$$ C_1 = \left ( \begin{array}{cccc}
\mathcal{R}_1^T  M \mathcal{R}_1 & & & \\
& \mathcal{R}_2^T M \mathcal{R}_2 & & \\
& & \ddots & \\
& & & \mathcal{R}_p^T A \mathcal{R}_p \end{array} \right ), \quad 
B_1 = \left ( \begin{array}{cccc}
\mathcal{R}_1^T  A \mathcal{R}_1 & & & \\
& \mathcal{R}_2^T A \mathcal{R}_2 & & \\
& & \ddots & \\
& & & \mathcal{R}_p^T A \mathcal{R}_p \end{array} \right ),$$
where $\{ \mathcal{R}_i \}_{i=1}^p$ satisfies 
$$ \mathcal{R}_{\text{off}} = \left ( \mathcal{R}_1 \ \mathcal{R}_2 \ \cdots \ \mathcal{R}_p \right ).$$


\section{Numerical results} \label{sec:numerics}
In this section, we present some numerical examples to illustrate the efficiency of the proposed splitting algorithm. 
In the numerical experiment, we set $\Omega = (0,1)^2$ and $T = 0.25$. Let $\mathcal{T}^H$ be a uniform coarse rectangular partition such that $\Omega$ is equally divided into $ 16 \times 16 $ small pieces of square with mesh size $H = \sqrt{2}/16$. We further divide each coarse-grid block $K \in \mathcal{T}^H$ into $16 \times 16$ uniform small squares. Each degree of freedom defined in coarse scale can be downscaled to a fine (uniform and rectangular) mesh $\mathcal{T}^h$ with mesh size $h = \sqrt{2}/256$. 

In the numerical experiments below, we take $C_1$ and $B_1$ to be the block-diagonal parts of the complete matrices $C$ and $B$, respectively. That is, we take $C_2 = C - C_1$ and $B_2 = B - B_1$ with 
$$ C_1 = \left ( \begin{array}{cc}
\mathcal{R}_1^T M \mathcal{R}_1 & O \\
O & \mathcal{R}_2^T M \mathcal{R}_2 
\end{array} \right ) \quad \text{and} \quad 
B_1 = \left ( \begin{array}{cc}
\mathcal{R}_1^T A \mathcal{R}_1 & O \\
O & \mathcal{R}_2^T A \mathcal{R}_2 
\end{array} \right ).$$

We remark that during the evolution of the solution process, one needs to solve a linear system of the following form for $n \in \mathbb{N}$ and $n \geq 2$: 
\begin{eqnarray}
\begin{split}
&(\mu C_1 + \tau \sigma B_1) Z^{n+1} =  \varphi^n,\\
&\varphi^n = \tau \tilde f^{n+1} - \Big [ \tau (1-\sigma) B_1 + \tau B_2 + (1-2\mu)C_1 + C_2 \Big ] Z^n + \Big [  (1-\mu) C_1 + C_2 \Big ] Z^{n-1}.
\end{split}
\label{eqn:com-reduced}
\end{eqnarray}

The left-hand side matrices $C_1$ and $B_1$ in the first equation above are block-diagonal. 
\revii{We briefly discuss the complexity reduction due to the proposed temporal splitting. In practical implementation, the solution $Z^{n+1}$ of equation \eqref{eqn:com-reduced} is obtained by solving two sub-problems with the left-hand sides $\mu C_{ii} + \tau \sigma B_{ii}$ (with $i = 1,2$) and concatenating the sub-solutions $V^{n+1}$ and $W^{n+1}$. In each iteration, the size of the linear system is reduced by half comparing the backward Euler scheme and one only needs to solve two smaller linear systems. As a result, the overall computational cost is reduced by half in each iteration. Also, these two sub-problems related to the left-hand side $\mu C_{ii} + \tau \sigma B_{ii}$ (with $i = 1,2$) are completely decoupled so one can achieve more computational saving with parallelization.}

To quantify the efficiency of the proposed solution decomposition algorithm, we record the errors between the multiscale solution $u_{\text{ms}}$ solved by the GMsFEM using the backward Euler method without any splitting and the one $U_{\text{split}}$ solved by the solution splitting algorithm. In particular, we define $$e_{L^2} := 
\frac{
\norm{(u_{\text{ms}} - U_{\text{split}})(T)}_{L^2(\Omega)}}
{\norm{u_{\text{ms}}(T)}_{L^2(\Omega)}}
\quad \text{and} \quad e_a := \frac{\norm{(u_{\text{ms}} - U_{\text{split}})(T)}_{a}}{\norm{u_{\text{ms}}(T)}_{a}}$$
to measure the efficiency of the proposed algorithm. Here, we define $\revi{\norm{v}_{L^2(\Omega)} := \left ( v,v \right )_\Omega^{1/2}}$ to be the $L^2$ norm for any $v \in L^2(\Omega)$ and $\norm{v}_a := \left ( \int_\Omega \kappa \abs{\nabla v}^2 ~ dx \right )^{1/2}$ the energy norm for any $v \in \revi{H_0^1(\Omega)}$. 

Furthermore, in all the experiments below, we will pick $\ell \in \mathbb{N}^+$ eigenfunctions (related to the first $\ell = \ell_i$ smallest eigenvalues obtained in \eqref{eqn:sp}) to form the local multiscale space corresponding to each coarse neighborhood $\omega_i$. 
To enhance the performance of computation, we perform a modified Gram-Schmidt process with respect to the energy $a(\cdot,\cdot)$ for these $\ell$ local multiscale basis functions to make them orthonormal. Under this setting, the matrix blocks $B_{12}$ and $B_{21}$ are highly sparse and therefore the whole splitting algorithm becomes more stable. 
In order to form $\mathcal{R}_1$ and $\mathcal{R}_2$ in the solution decomposition algorithm, we choose the first $\ell^{(1)} \in \mathbb{N}^+$ (with $\ell^{(1)} < \ell$) eigenfunction(s) corresponding to the first $\ell^{(1)}$ smallest eigenvalue(s) in each local multiscale space and form the matrix $\mathcal{R}_1$. 
Then, the rest of $\ell^{(2)} := \ell - \ell^{(1)}$ eigenfunction(s) in each local multiscale space is/are used to form another matrix $\mathcal{R}_2$. We denote $\ell^{(1)} + \ell^{(2)}$ such solution decomposition setting. 

\begin{example}\label{exp1}
In the first example, we consider a periodic permeability $\kappa$ to be as follows:
$$ \kappa(x_1, x_2) = \frac{2+ \sin(11\pi x_1) \sin(13 \pi x_2)}{1.4 + \cos (12 \pi x_1) \cos (7 \pi x_2)}$$
for any $(x_1, x_2) \in \Omega$. We set the source function to be 
$ \mathcal{F}(x_1, x_2) = \exp((x_1-0.5)^2+(x_2-0.5)^2)$
for any $(x_1, x_2) \in \Omega$ and the initial condition is $U_0(x_1, x_2) = \sin(\pi x_1)\sin(\pi x_2)$ for any $(x_1, x_2) \in \Omega$. We set $\tau = 10^{-3}$. 
In each local offline space of GMsFEM, we pick $\ell = 6$ eigenfunctions (related to the first $6$ smallest eigenvalues) to form the local multiscale space. 
In this case, the total degrees of freedom is $\mathcal{N}_{\text{c}} = 1350$. 
The parameters $(\mu ,\sigma) =(1,1)$ are chosen. 
The record of errors at terminal time with different splitting options is presented in Table \ref{tab:error-exp1} and 
the profiles of the solutions are sketched in Figure \ref{fig:sol-exp1}. 
The error at terminal time is about the magnitude of $10^{-4}$ in terms of both the $L^2$ and energy errors. 
One can observe that the proposed solution decomposition method works well in this simple case with periodic heterogeneous medium. 

\begin{figure}[ht!]
\centering
\includegraphics[width=2.5in]{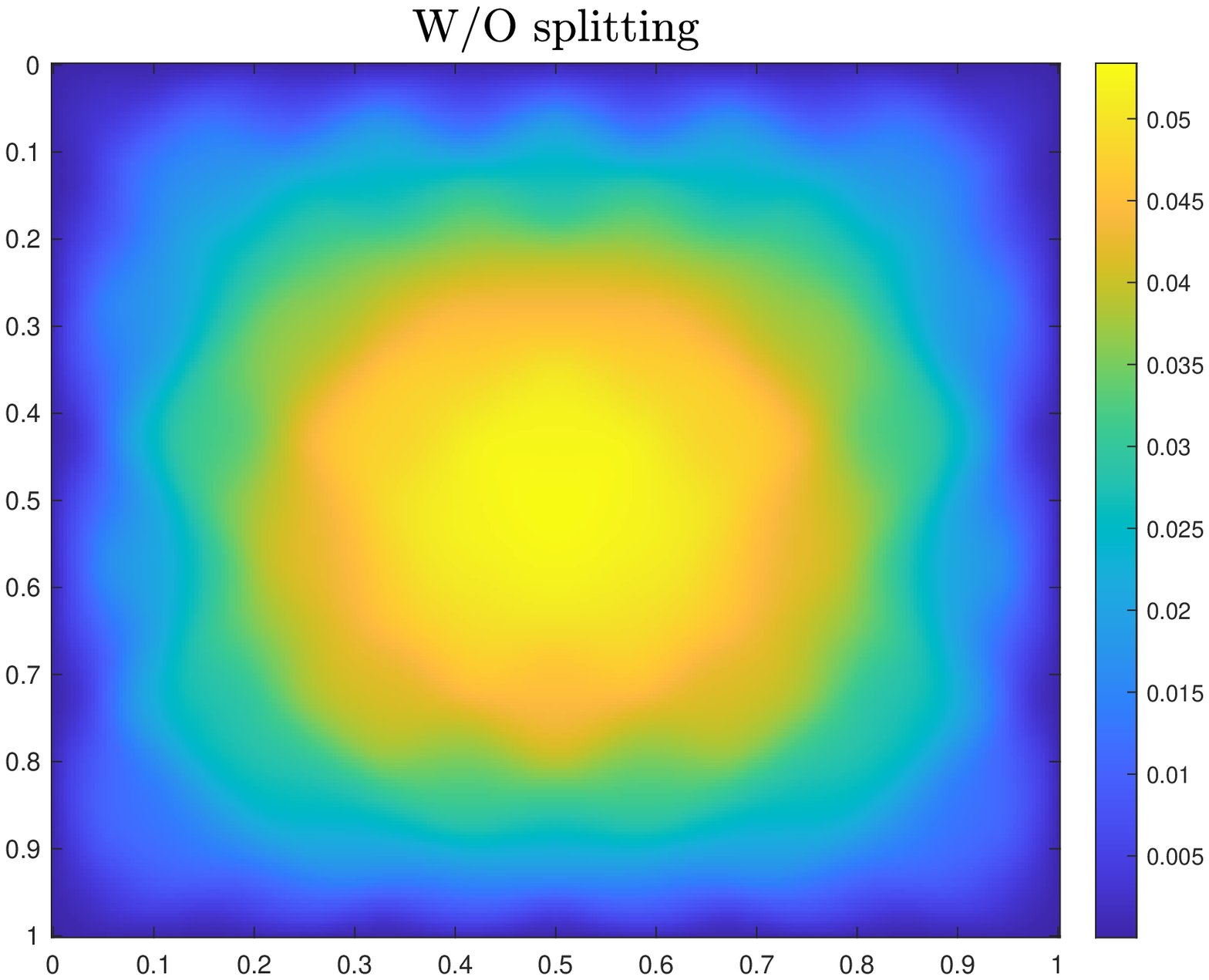}
\qquad 
\includegraphics[width=2.5in]{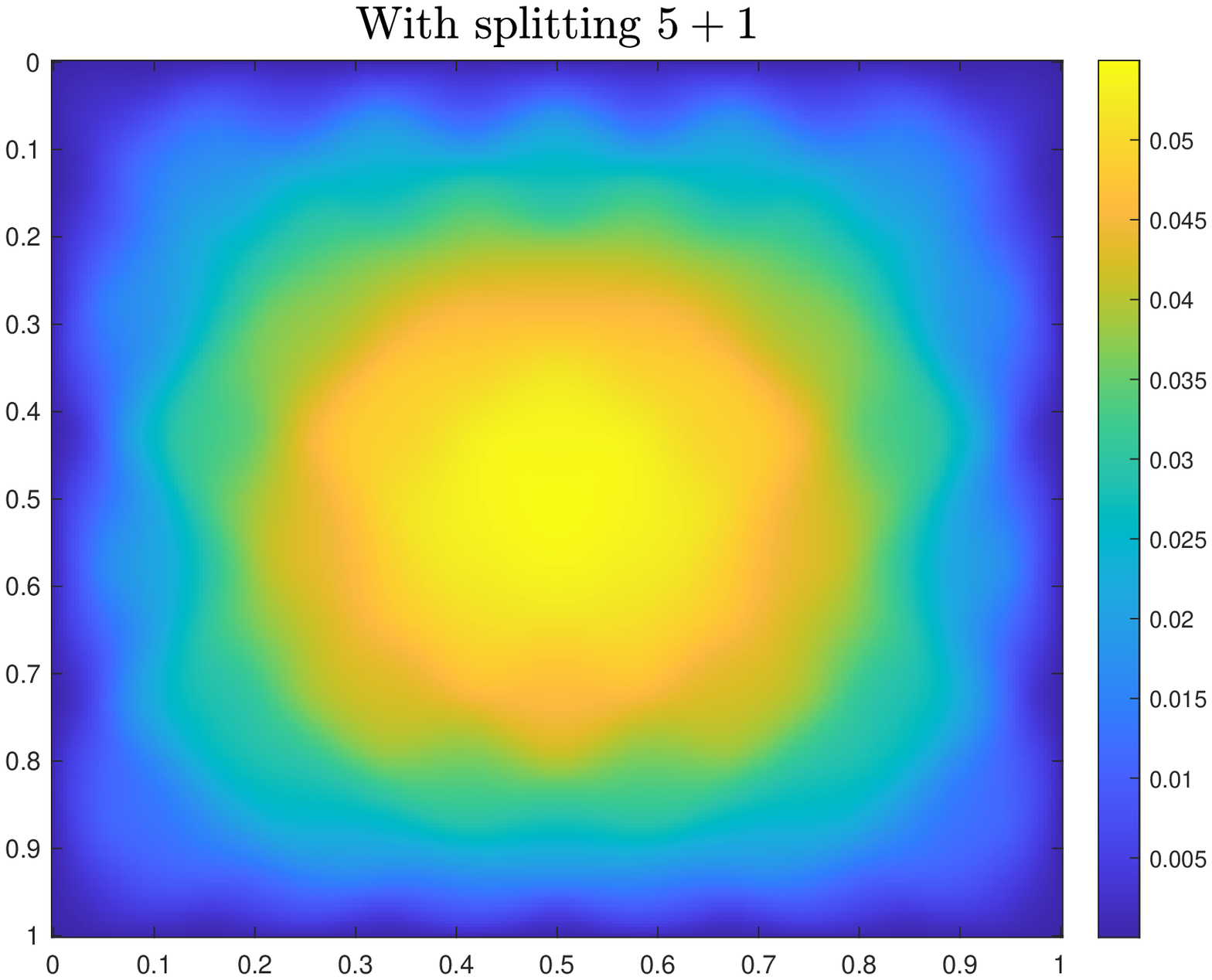}
\caption{Solutions at $T=0.25$. Left: no splitting. Right: with $5+1$ splitting (Example \ref{exp1}).} 
\label{fig:sol-exp1}
\end{figure}

\begin{table}[ht!]
\centering
\begin{tabular}{c|c|c}
\hline \hline 
$\ell^{(1)} + \ell^{(2)}$ & $e_{L^2}$ & $e_a$ \\
\hline \hline 
$1+5$ & $6.3740\times 10^{-5}$ & $6.1825\times 10^{-5}$ \\ 
$2+4$ & $5.0456\times 10^{-4}$ & $4.7575\times 10^{-4}$ \\ 
$3+3$ & $4.8084\times 10^{-4}$ & $5.7561\times 10^{-4}$ \\ 
$4+2$ & $4.8090\times 10^{-4}$ & $7.3801\times 10^{-4}$ \\ 
$5+1$ & $4.7101\times 10^{-4}$ & $4.4519\times 10^{-4}$ \\ 
\hline \hline 
\end{tabular}
\caption{Errors at $T=0.25$ with $H = \frac{\sqrt{2}}{16}$, $\tau = 10^{-3}$, and $(\mu, \sigma) = (1,1)$ (Example \ref{exp1}).}
\label{tab:error-exp1}
\end{table}

\begin{figure}[ht!]
\centering
\includegraphics[width=2.7in]{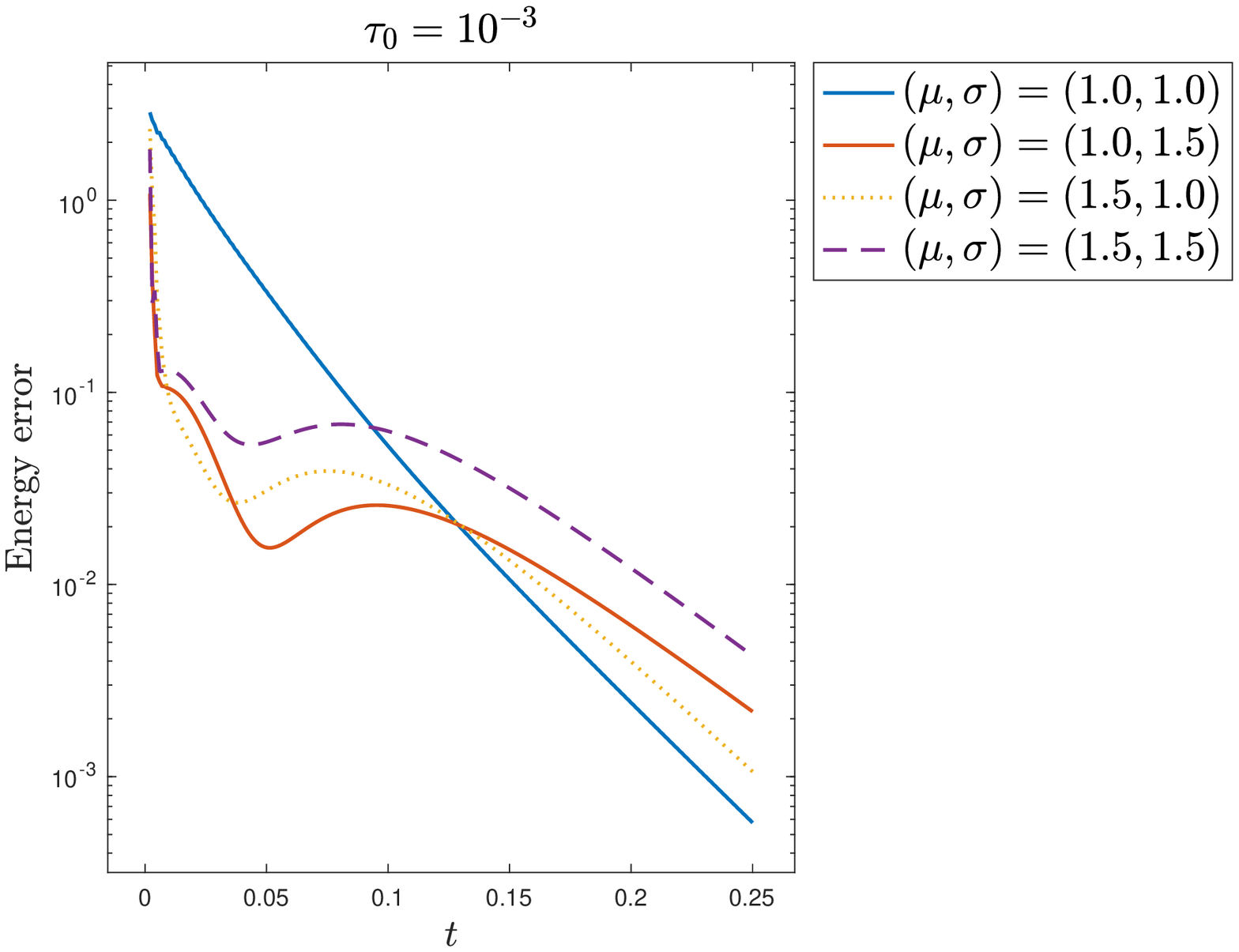}
\qquad 
\includegraphics[width=2.7in]{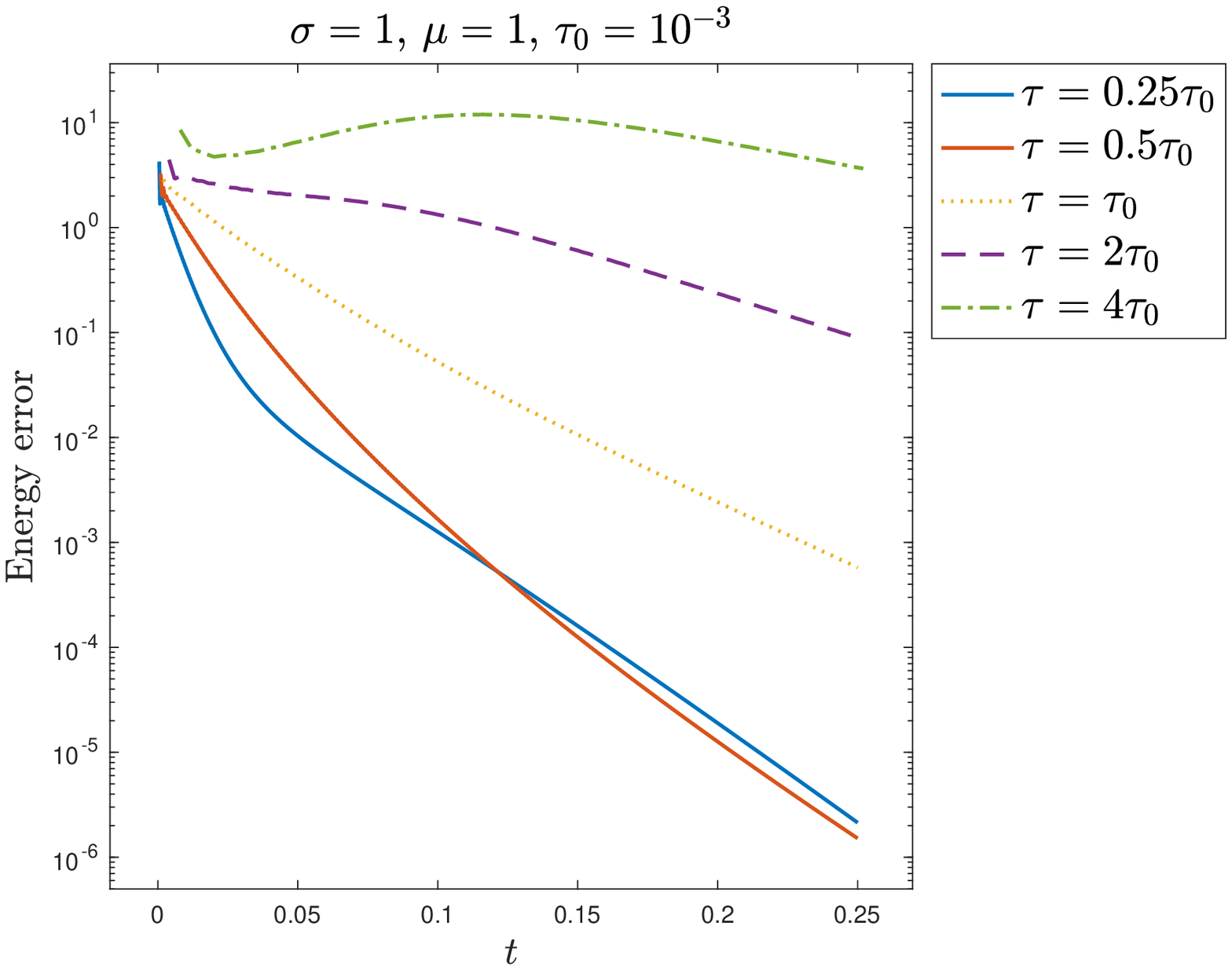}
\caption{History of energy error (using $3+3$ splitting). Left: $\sigma  \in \{ 1.0, 1.5 \}$ and $\mu \in \{ 1.0, 1.5\}$ with $\tau_0= 10^{-3}$. Right: $\tau \in \{ 0.25 \tau_0, 0.5 \tau_0, \tau_0, 2\tau_0, 4\tau_0 \}$ with $(\mu, \sigma ) = (1,1)$ (Example \ref{exp1}). }
\label{fig:his-exp1}
\end{figure}

We explore the dependence of the parameters $\mu$ and $\sigma$ for the solution decomposition algorithm. To illustrate this, we report the history of energy error against time evolution using a $3+3$ splitting. In Figure \ref{fig:his-exp1} (left), the history of energy error is sketched for four different settings of the parameters $\mu$ and $\sigma$. The splitting scheme with $(\mu, \sigma) = (1,1)$ performs more stable in the error reduction against time evolution. 
In Figure \ref{fig:his-exp1} (right), we fix $(\mu, \sigma )= (1,1)$ and set different temporal step size $\tau \in \{ 0.25 \tau_0, 0.5 \tau_0, \tau_0, 2\tau_0, 4\tau_0 \}$ with $\tau_0 = 10^{-3}$. The convergence with respect to the temporal step size can be observed. 
\end{example}

\begin{example} \label{exp2}
In this example, we consider a more heterogeneous permeability with high-contrast features, which is more challenging. The permeability is depicted in Figure \ref{fig:exp2} (left). The (time-independent) source function is sketched in Figure \ref{fig:exp2} (right). The initial condition is the same as in Example \ref{exp1} and we set $\tau = 2\times 10^{-4}$. 
We take $\ell = 10$ eigenfunctions corresponding the first $10$ smallest eigenvalues to form the local multiscale space. In this setting, the degrees of freedom is $\mathcal{N}_{\text{c}} = 2250$. The parameters of splitting $(\mu, \sigma) = (1, 1)$ are chosen. 
The record of errors is presented in Table \ref{tab:error-exp2}. 
The profiles of the multiscale solution without splitting and the one with $6+4$ splitting are sketched in Figure \ref{fig:sol-exp2}. 

\begin{figure}[ht!]
\centering
\includegraphics[width=2.5in]{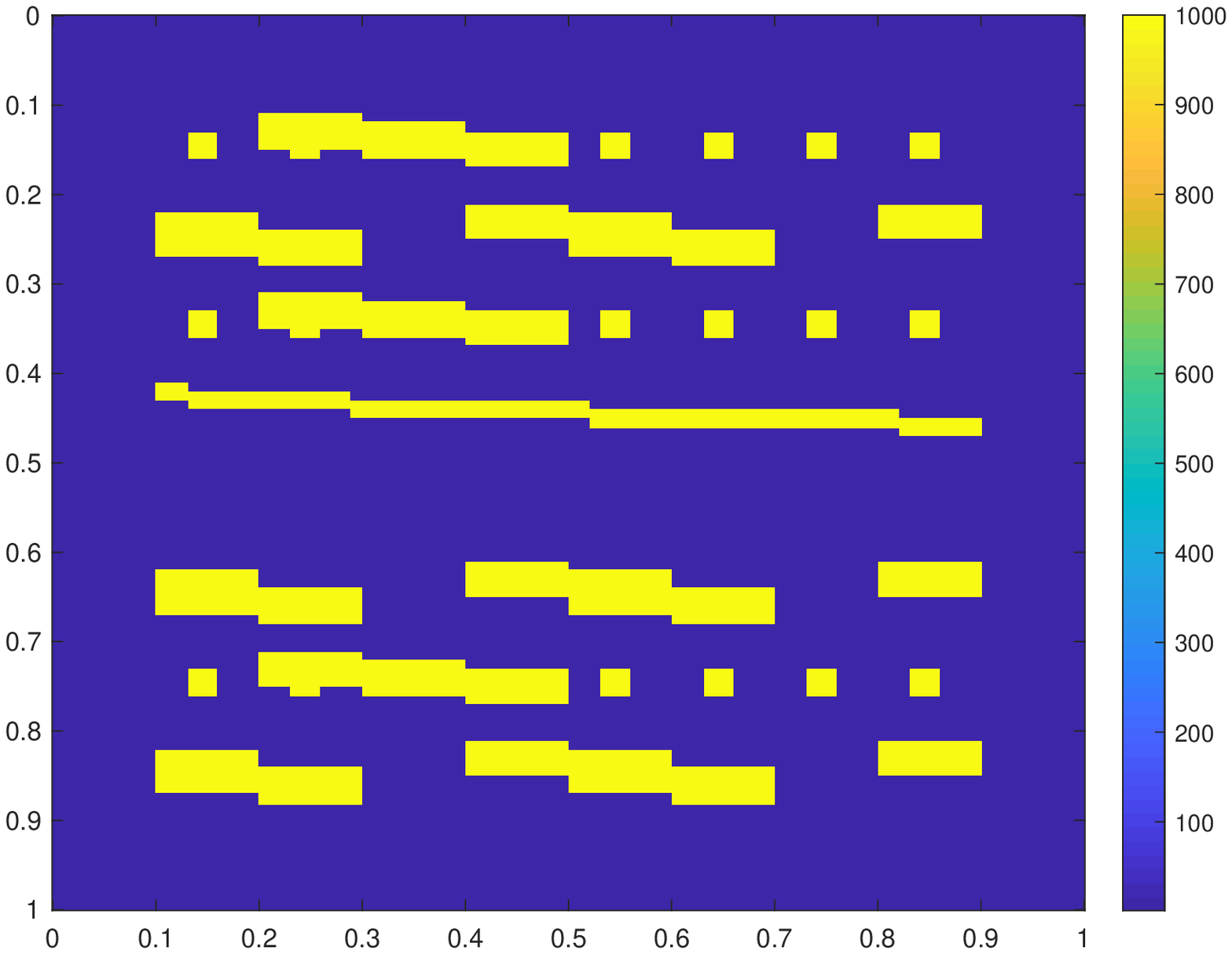}
\qquad 
\includegraphics[width=2.5in]{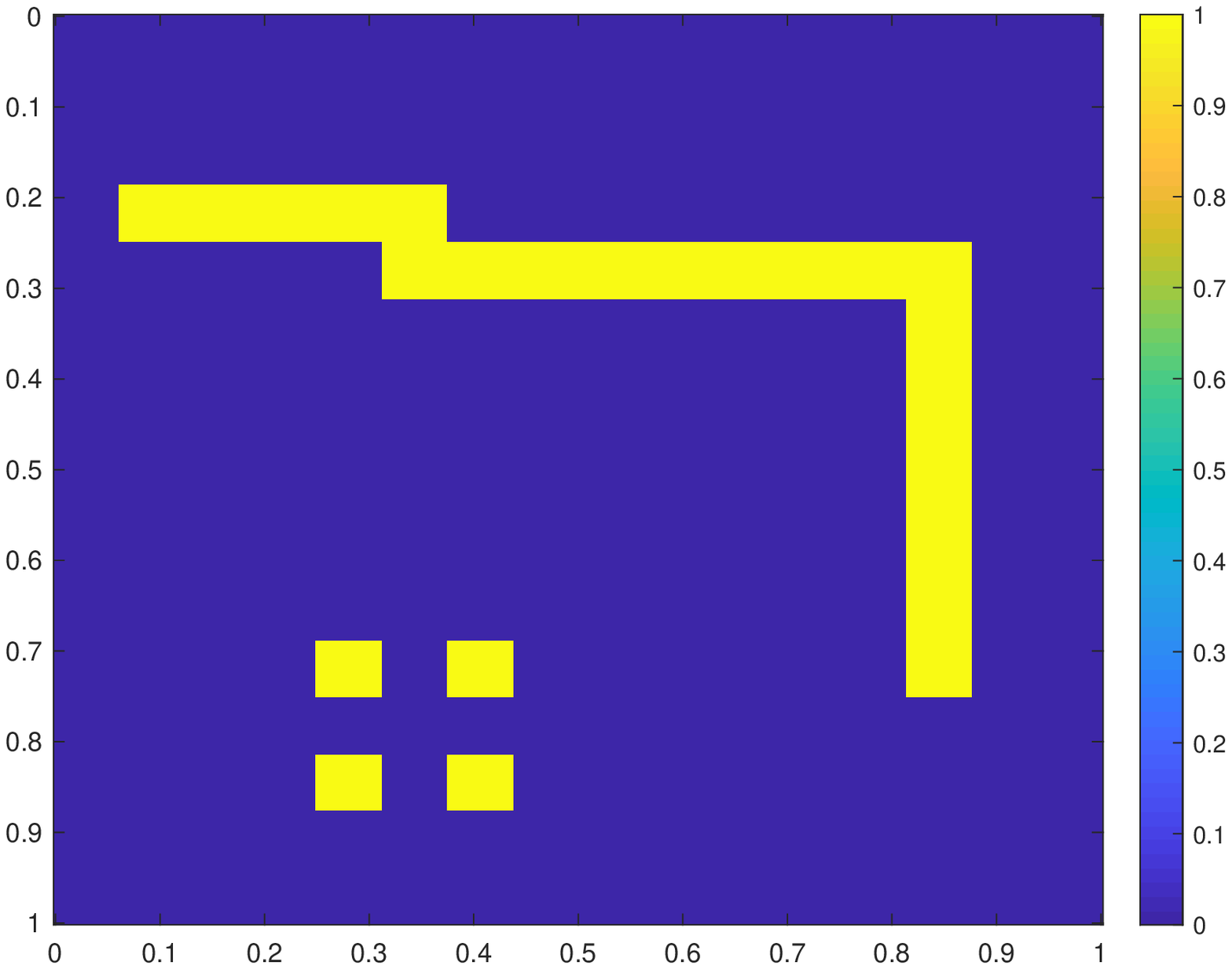}
\caption{Problem setting in Example \ref{exp2}. Left: the permeability. Right: the source function.} 
\label{fig:exp2}
\end{figure}

From Table \ref{tab:error-exp2}, we observe that the splitting $1+9$ outperforms the rest of the options of splitting. 
We remark that the splitting algorithm fails to converge if $\mu < 1.0$ and $\sigma < 1.0$ in this example. 
Moreover, the $8+2$ and $9+1$ splittings are moderately time-consuming than the rest of the splitting schemes, since the main matrices $C_1$ and $B_1$ are of the largest band-width among all the options of splitting. 

\begin{figure}[ht!]
\centering
\includegraphics[width=2.5in]{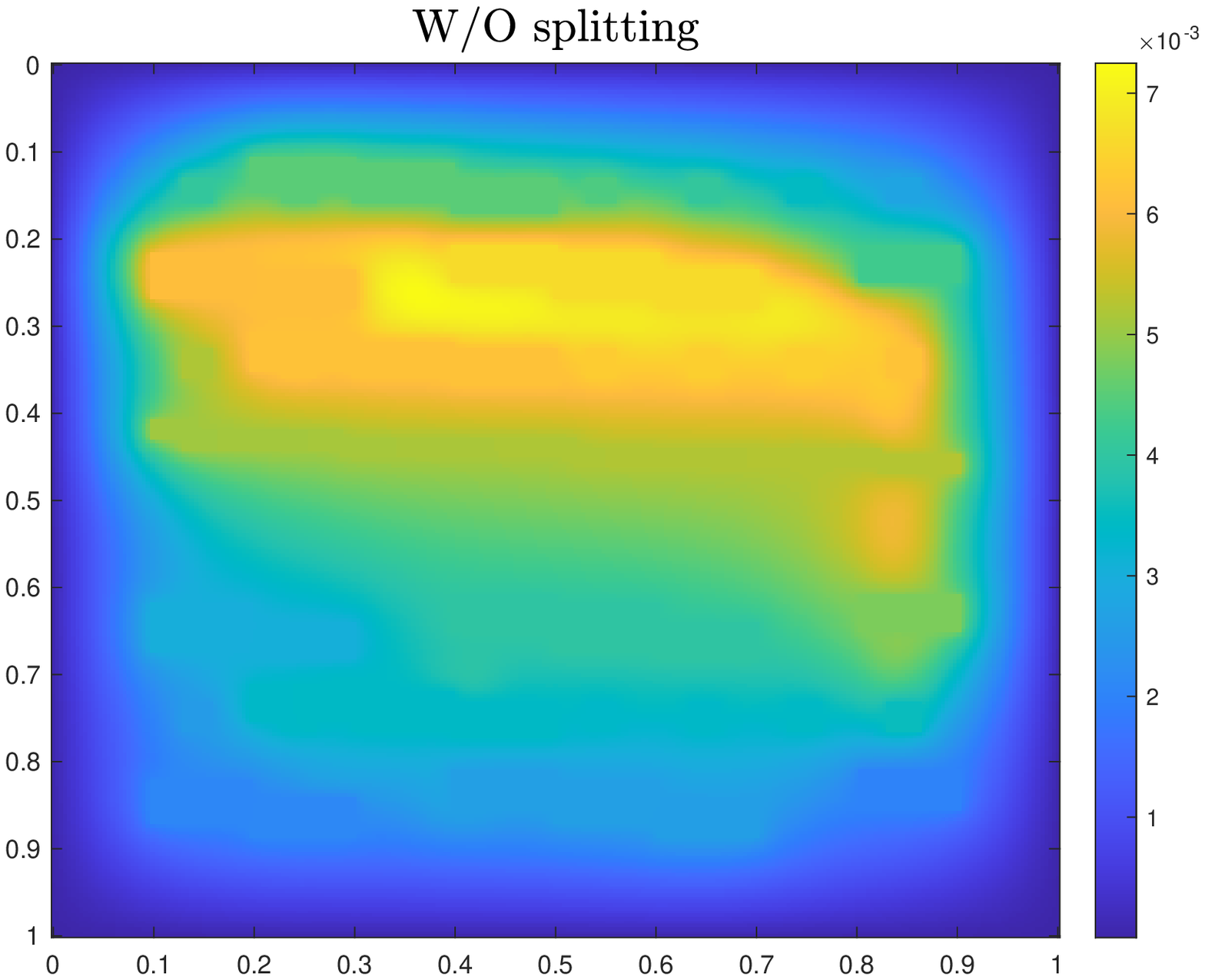}
\qquad 
\includegraphics[width=2.5in]{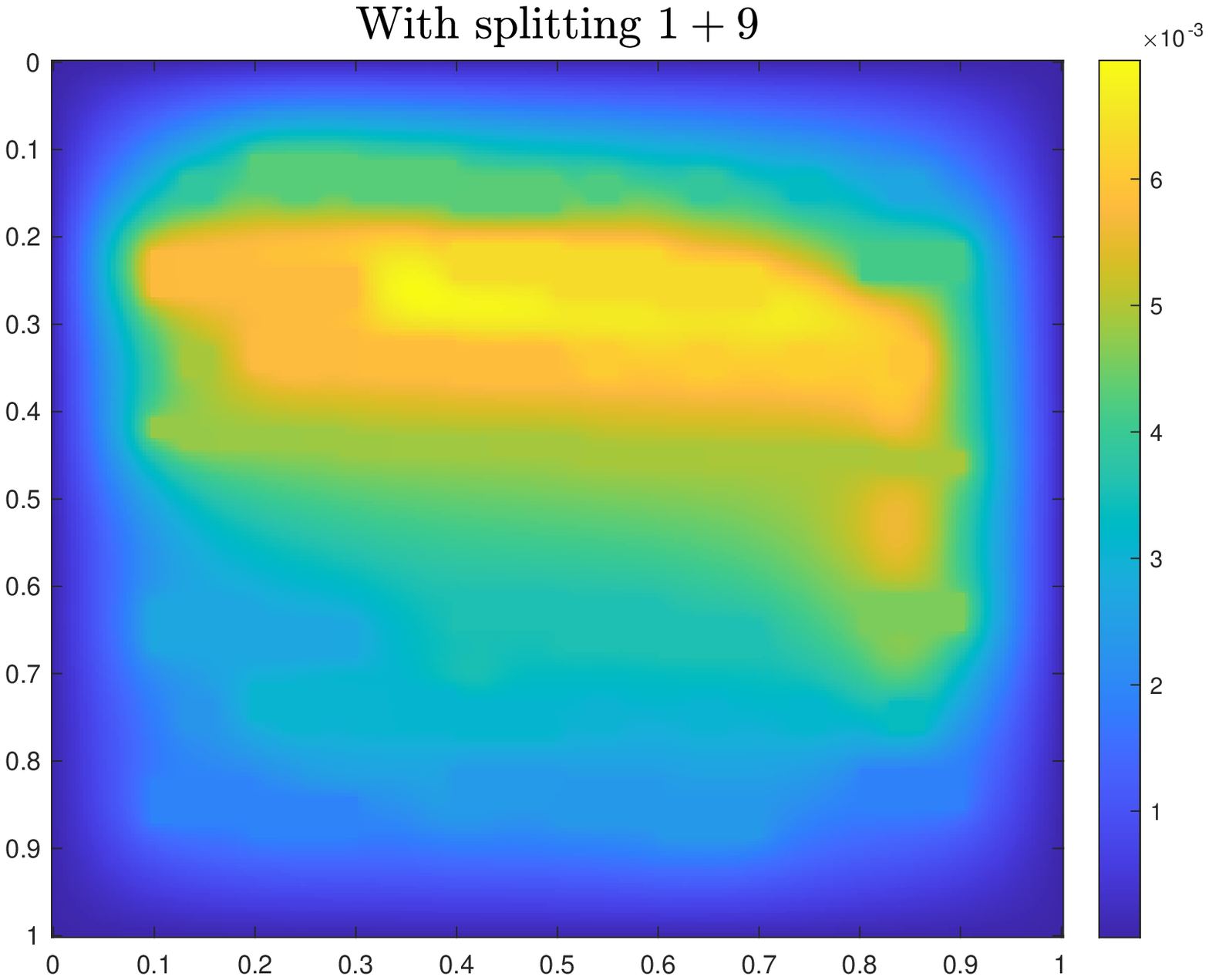}
\caption{Solution at $T = 0.25$. Left: no splitting. Right: with $1+9$ splitting (Example \ref{exp2}).} 
\label{fig:sol-exp2}
\end{figure}

\begin{table}[ht!]
\centering
\begin{tabular}{c|c|c}
\hline \hline 
$\ell^{(1)} + \ell^{(2)}$ & $e_{L^2}$ & $e_a$ \\
\hline \hline 
$1+9$ & $3.9530\times 10^{-3}$ & $3.5080\times 10^{-3}$ \\ 
$2+8$ & $5.7982\times 10^{-2}$ & $5.1215\times 10^{-2}$ \\ 
$3+7$ & $0.1337$ & $0.1175$ \\ 
$4+6$ & $0.2587$ & $0.2323$ \\ 
$5+5$ & $0.3819$ & $0.3483$ \\ 
$6+4$ & $0.4386$ & $0.3968$ \\ 
$7+3$ & $0.4510$ & $0.4315$ \\ 
$8+2$ & $0.4580$ & $0.4672$ \\ 
$9+1$ & $0.3517$ & $0.3257$ \\ 
\hline \hline 
\end{tabular}
\caption{Errors at $T=0.25$ with $H = \frac{\sqrt{2}}{16}$, $\tau = 2 \times 10^{-4}$, and $(\mu,\sigma) = (1,1)$ (Example \ref{exp2}).}
\label{tab:error-exp2}
\end{table}

In Figure \ref{fig:his-exp2}, we report the history of energy error using the $1+9$ splitting. We compare different settings of the parameters in Figure \ref{fig:his-exp2} (left). Among all options of parameters, the one with $(\mu, \sigma) = (1.5,1.0)$ (the dotted line in Figure \ref{fig:his-exp2} (left)) outperforms others in terms of the energy error. 
We also test with different options of temporal step size.  We denote $\tau_0 = 2 \times 10^{-4}$ and we use $1+9$ splitting with $\tau \in \{ 0.25 \tau_0, 0.5 \tau_0, \tau_0, 2 \tau_0, 4\tau_0 \}$. The corresponding results are shown in Figure \ref{fig:his-exp2} (right). One can observe the convergence with respect to the temporal step size. 

\begin{figure}[ht!]
\centering
\includegraphics[width=2.7in]{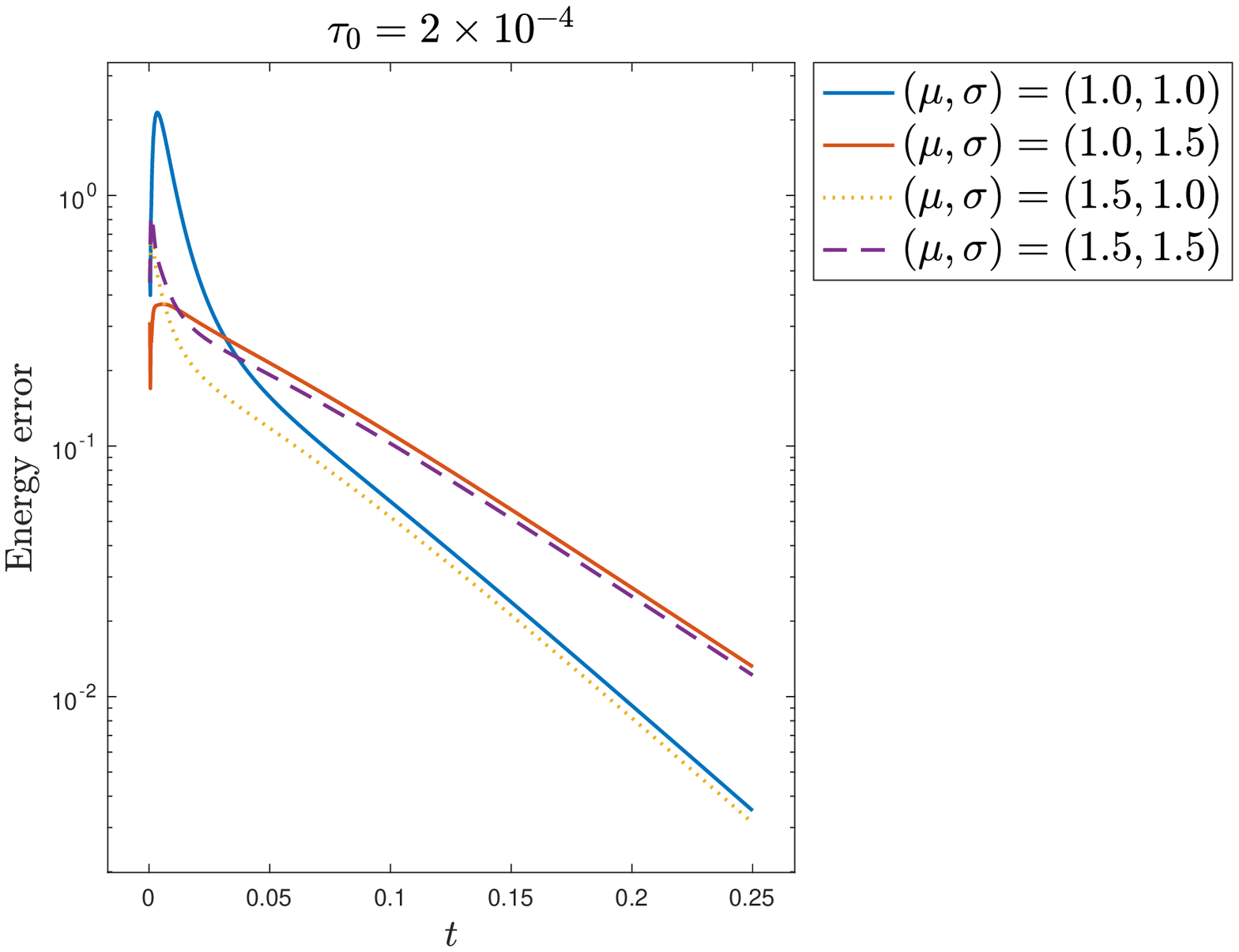}
\qquad 
\includegraphics[width=2.7in]{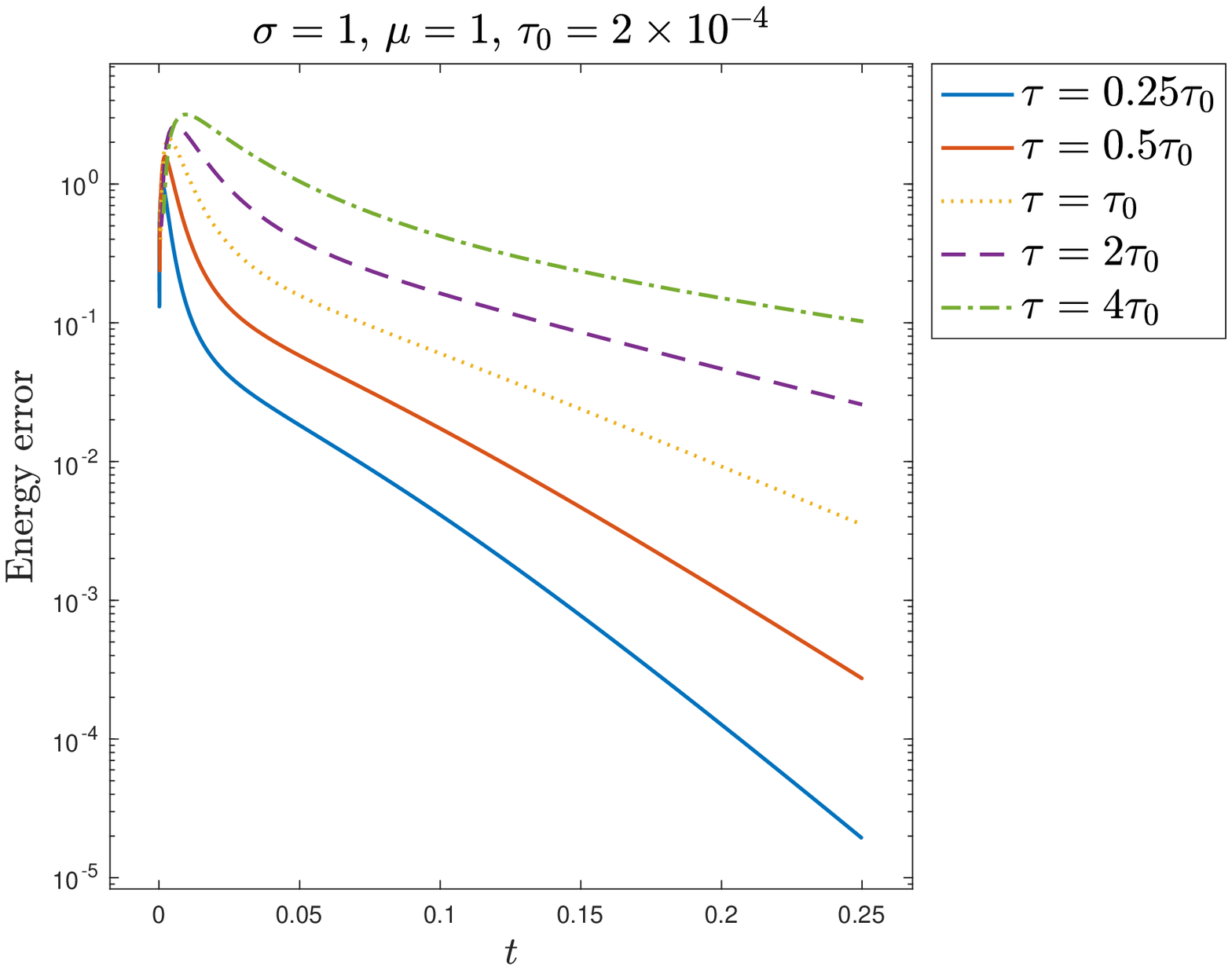}
\caption{History of energy error (using $1+9$ splitting). Left: $\sigma  \in \{ 1.0, 1.5 \}$ and $\mu \in \{ 1.0, 1.5\}$ with $\tau_0 = 2\times10^{-4}$. Right: $\tau \in \{ 0.25 \tau_0, 0.5 \tau_0, \tau_0, 2\tau_0, 4\tau_0 \}$ with $(\mu, \sigma ) = (1,1)$ (Example \ref{exp2}).
}
\label{fig:his-exp2}
\end{figure}
\end{example}

\begin{example}[Time-dependent source function] \label{exp3}
In this example, we consider a case with a time-dependent source function in the model problem \eqref{eqn:model-pde}. We set the source function to be 
$$ \mathcal{F}(t,x_1,x_2) = (\sin(\pi t)+1)\sin(\pi x_1)\sin(\pi x_2)$$
for any $t \in [0,T]$ and $(x_1, x_2) \in \Omega$. The permeability field considered in this example has long and thin channelized features in some part of the domain and it is depicted in Figure \ref{fig:setting-exp3} (left). The initial condition in this example is the same as in Examples \ref{exp1} and \ref{exp2}. We set $\tau = 10^{-3}$. 

\begin{figure}[ht!]
\centering
\includegraphics[width=2.5in]{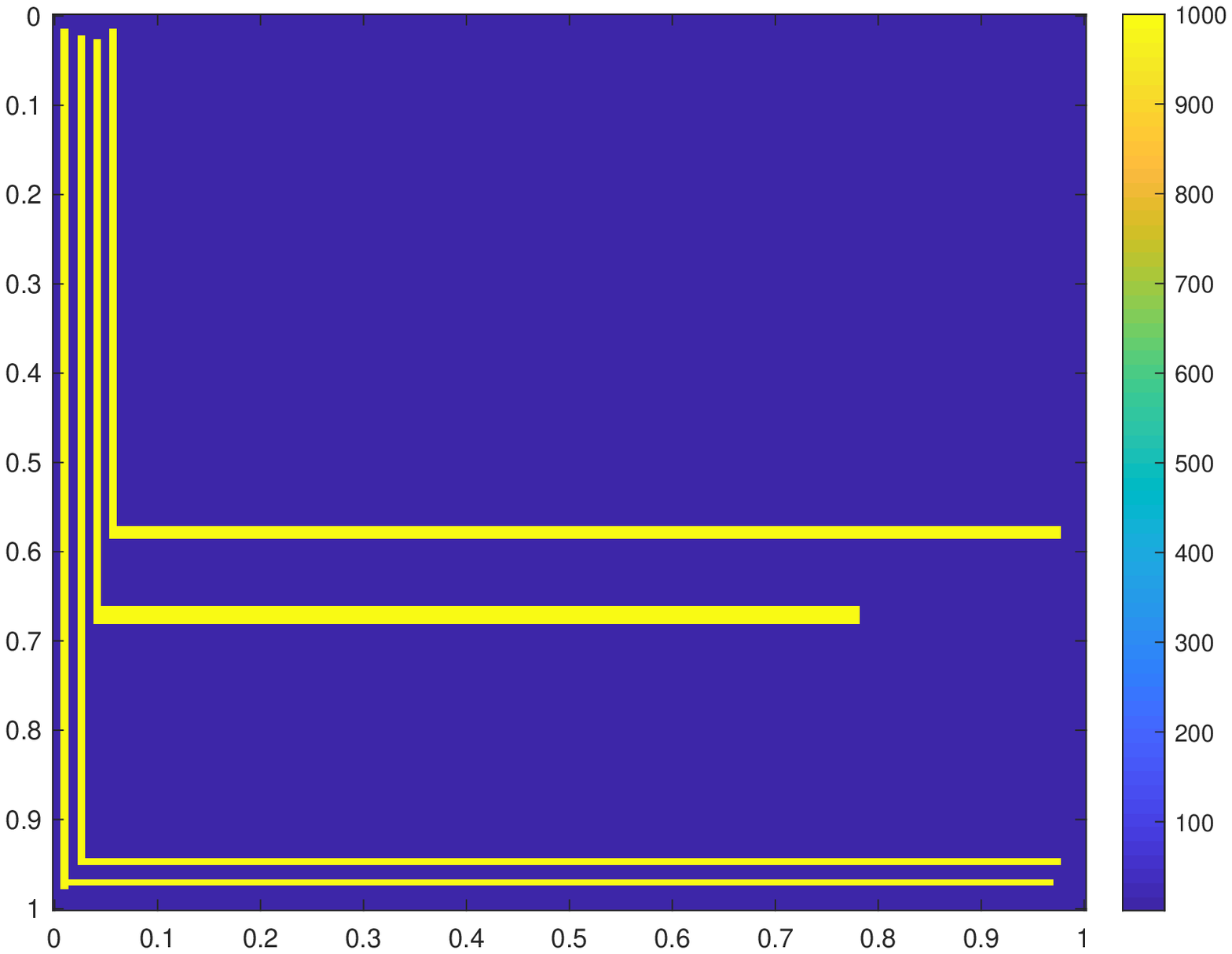}
\qquad 
\includegraphics[width=2.5in]{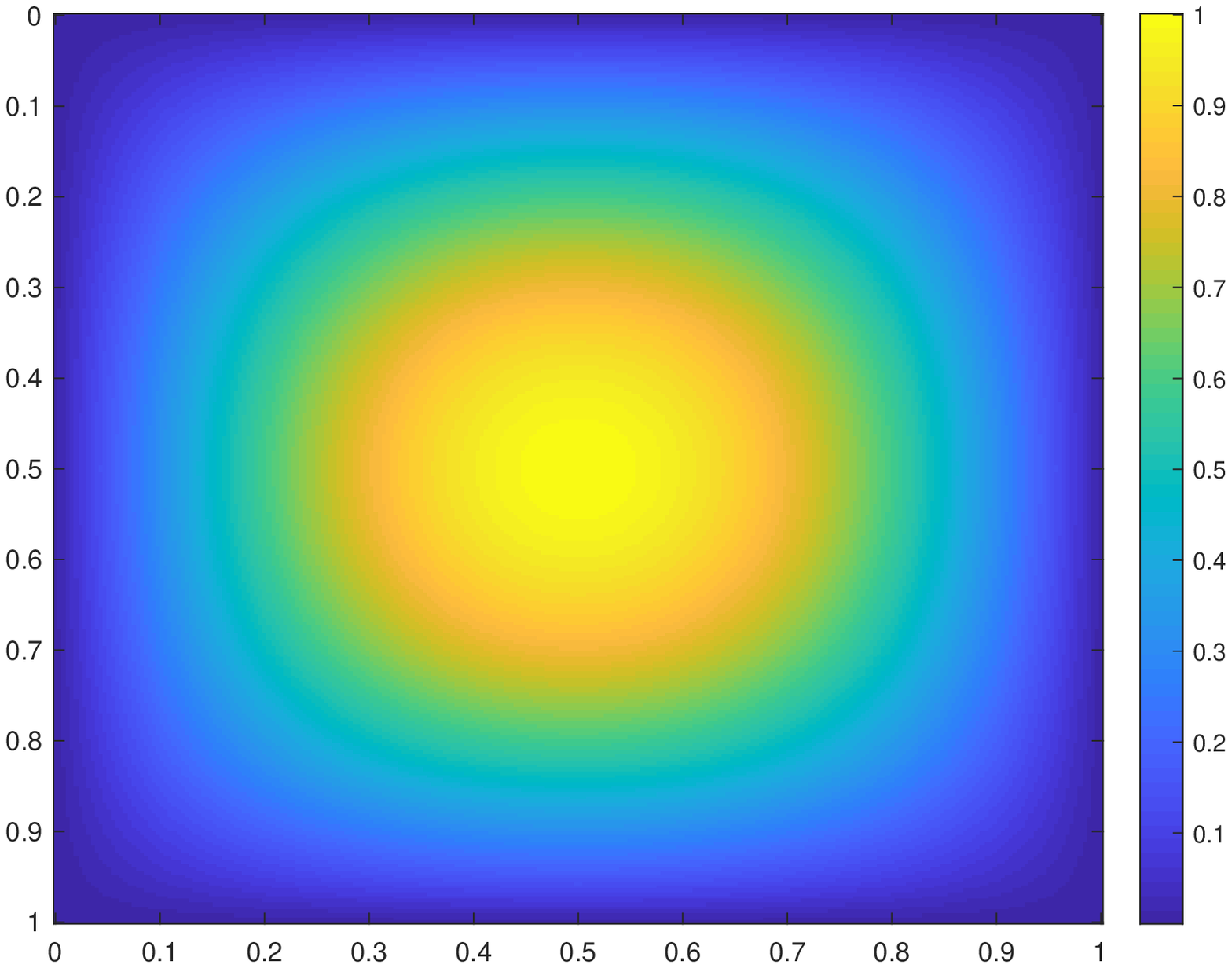}
\caption{Problem setting in Example \ref{exp3}. Left: the permeability. Right: the source function at $t = 0$.}
\label{fig:setting-exp3}
\end{figure}

Table \ref{tab:error-exp3} records the energy errors at terminal time with different splitting options. 
The profiles of the solution without splitting and the one with $6+4$ splitting are sketched in Figure \ref{fig:sol-exp3}. In this example, we observe that the $1+9$ splitting outperform others and the rest of different splitting options are around the same level of magnitude in terms of the errors at the terminal time. We remark that one may use a $\ell^{(1)} + \ell^{(2)}$ splitting with $\ell^{(1)} \ll \ell^{(2)}$ to speed up the computation process. 

\begin{table}[ht!]
\centering
\begin{tabular}{c|c|c}
\hline \hline 
$\ell^{(1)} + \ell^{(2)}$ & $e_{L^2}$ & $e_a$ \\
\hline \hline 
$1+9$ & $1.7779\times 10^{-6}$ & $3.9749\times 10^{-4}$ \\ 
$2+8$ & $2.3798\times 10^{-4}$ & $2.6205\times 10^{-3}$ \\ 
$3+7$ & $3.0263\times 10^{-4}$ & $2.0204\times 10^{-3}$ \\ 
$4+6$ & $4.6231\times 10^{-4}$ & $2.5952\times 10^{-3}$ \\ 
$5+5$ & $7.2235\times 10^{-4}$ & $3.4967\times 10^{-3}$ \\ 
$6+4$ & $7.3497\times 10^{-4}$ & $4.1353\times 10^{-3}$ \\ 
$7+3$ & $7.4922\times 10^{-4}$ & $6.2690\times 10^{-3}$ \\ 
$8+2$ & $1.0566\times 10^{-3}$ & $3.3773\times 10^{-2}$ \\ 
$9+1$ & $1.1311\times 10^{-3}$ & $1.9913\times 10^{-2}$ \\ 
\hline \hline 
\end{tabular}
\caption{Errors at $T=0.25$ with $H = \frac{\sqrt{2}}{16}$, $\tau = 2.5\times 10^{-4}$, and $(\mu,\sigma) = (1,1)$ (Example \ref{exp3}).}
\label{tab:error-exp3}
\end{table}

We also compare the performance using different parameters $(\mu,\sigma)$ for a fixed $1+9$ splitting option. 
The corresponding results are depicted in Figure \ref{fig:his-exp3} (left). 
In particular, it is noticeable that the case of $(\mu, \sigma) = (1,1)$ (i.e. the blue solid line in Figure \ref{fig:his-exp3} (left)) serves as a baseline of the performance of the splitting algorithm and the rest of splitting schemes work better than the baseline case. 
Denoting $p$ the number in the additive representation, the theoretical statement in \cite[Theorem 1]{efendiev2020splitting} suggests that one needs to set the parameters such that they are larger than $1$, which is $p/2$ in this situation with $p =2$. 

\begin{figure}[ht!]
\centering
\includegraphics[width=2.5in]{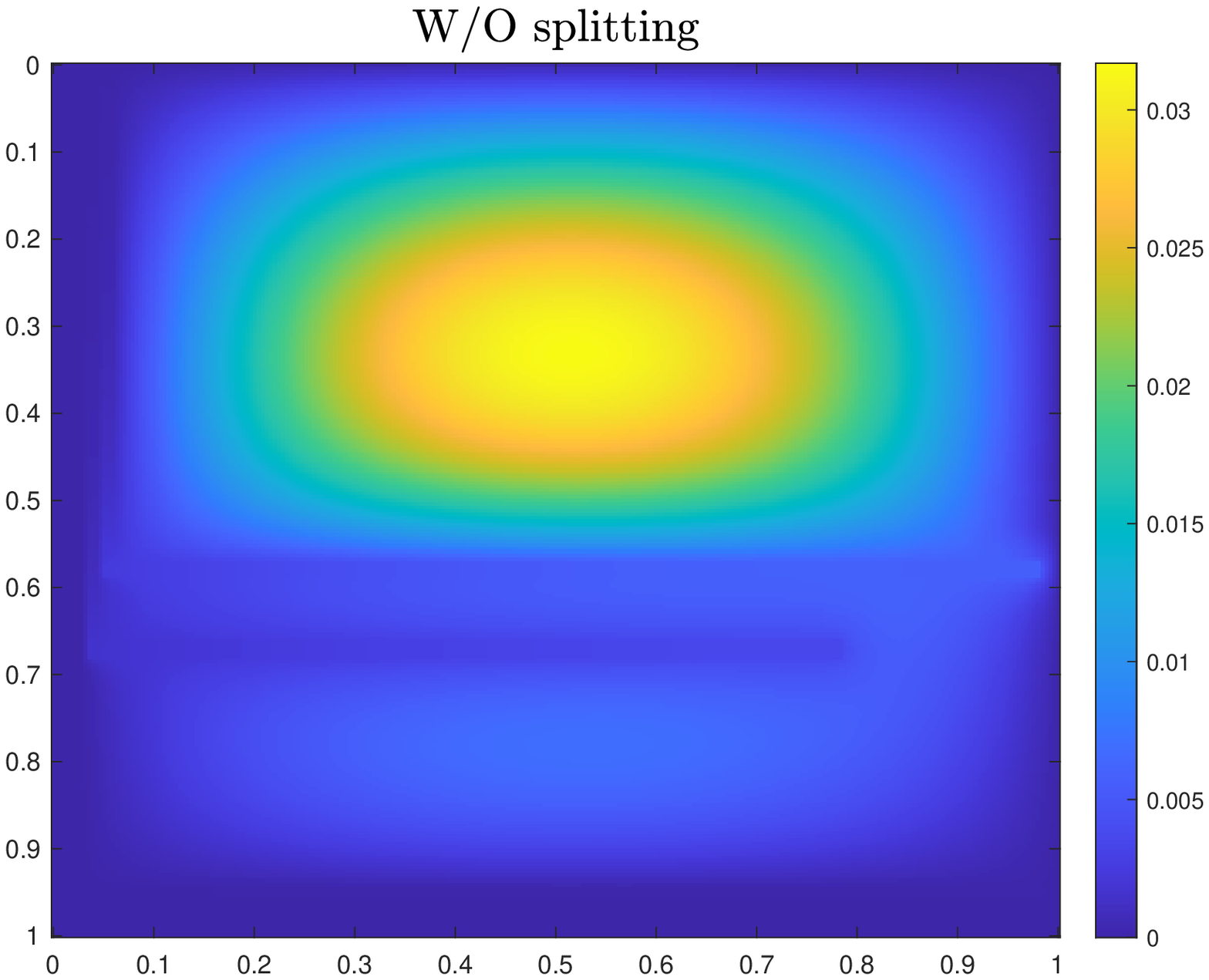}
\qquad 
\includegraphics[width=2.5in]{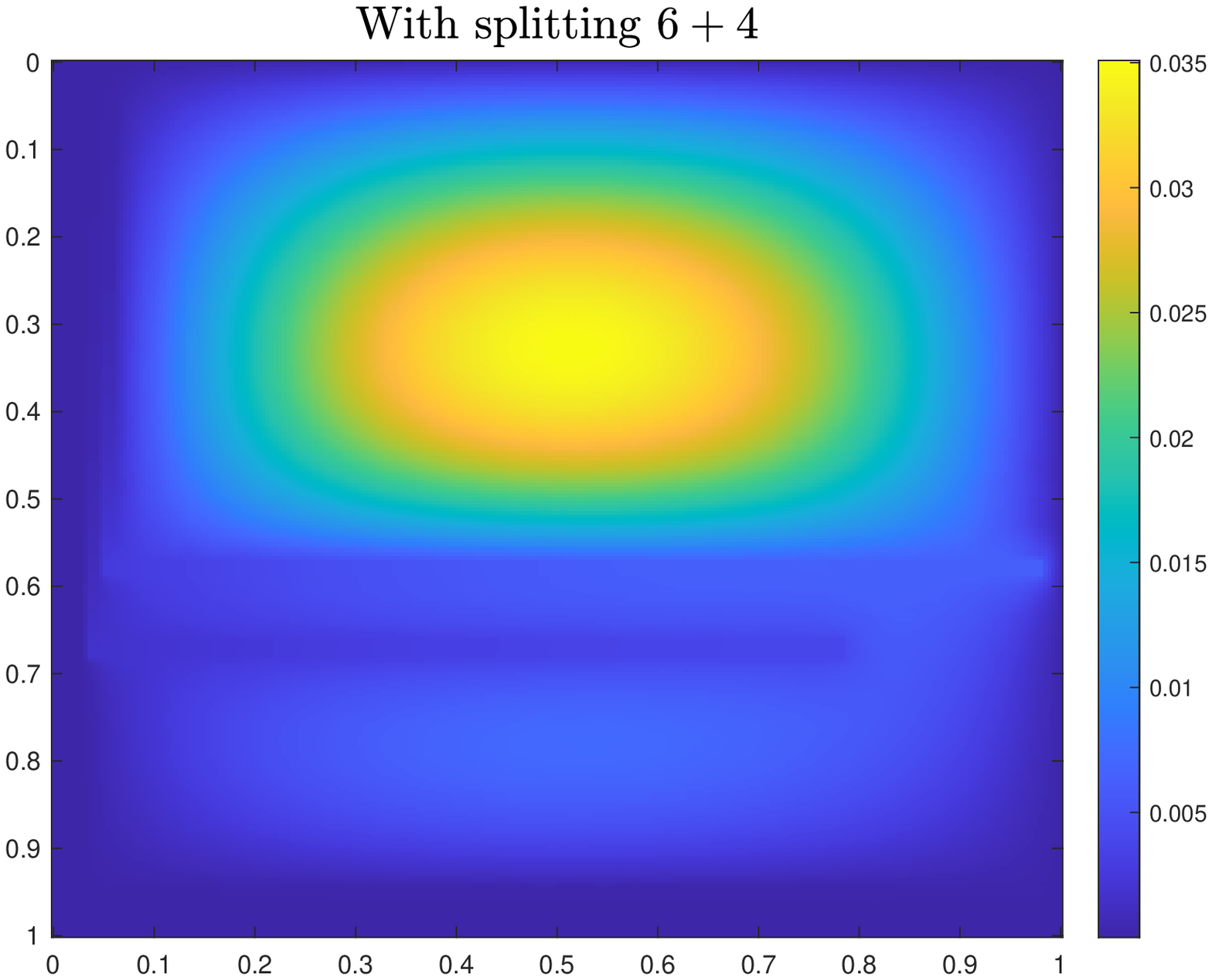}
\caption{Solution at $T=0.25$. Left: no splitting. Right: with $6+4$ splitting (Example \ref{exp3}).} 
\label{fig:sol-exp3}
\end{figure}

\begin{figure}[ht!]
\centering
\includegraphics[width=2.7in]{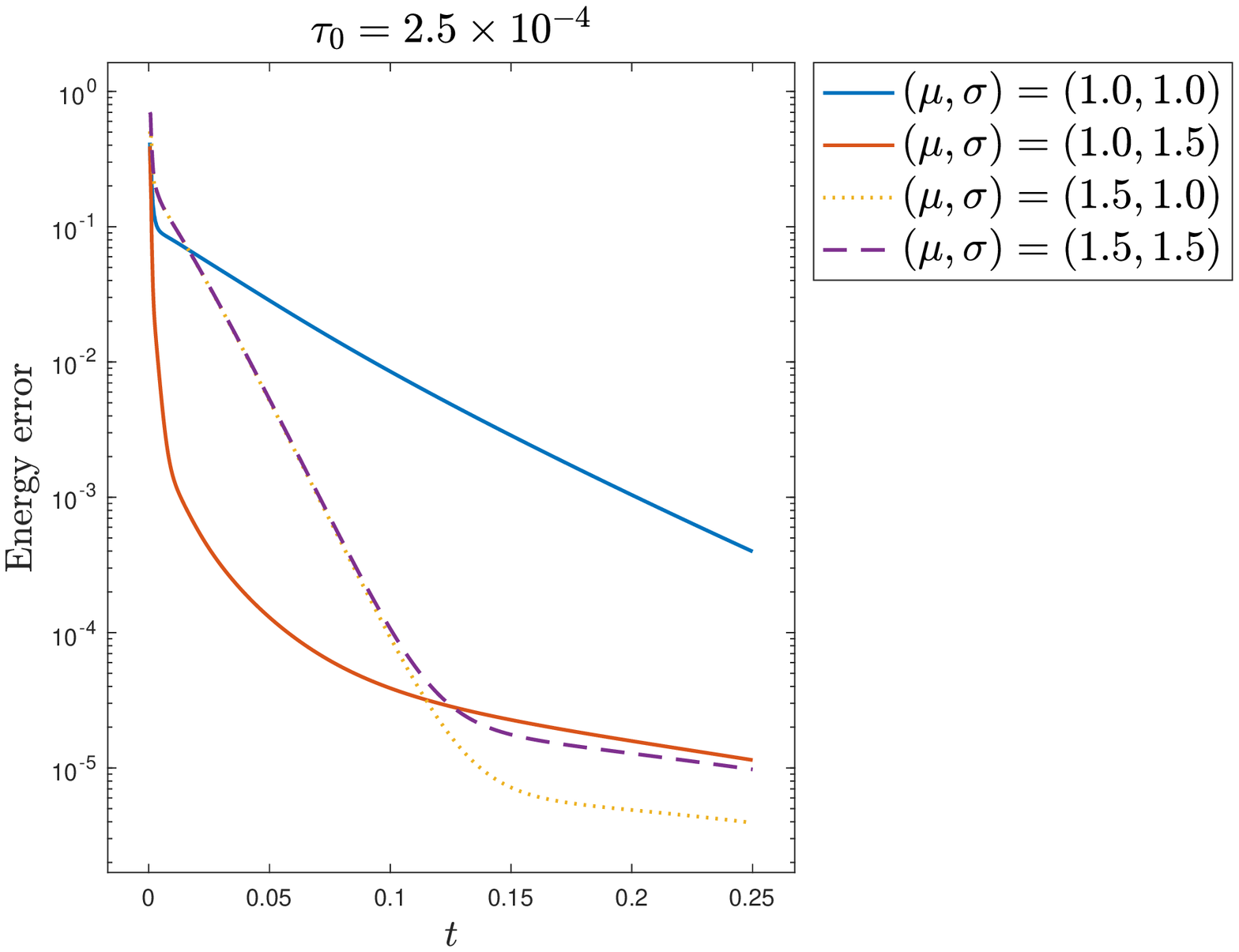}
\qquad 
\includegraphics[width=2.7in]{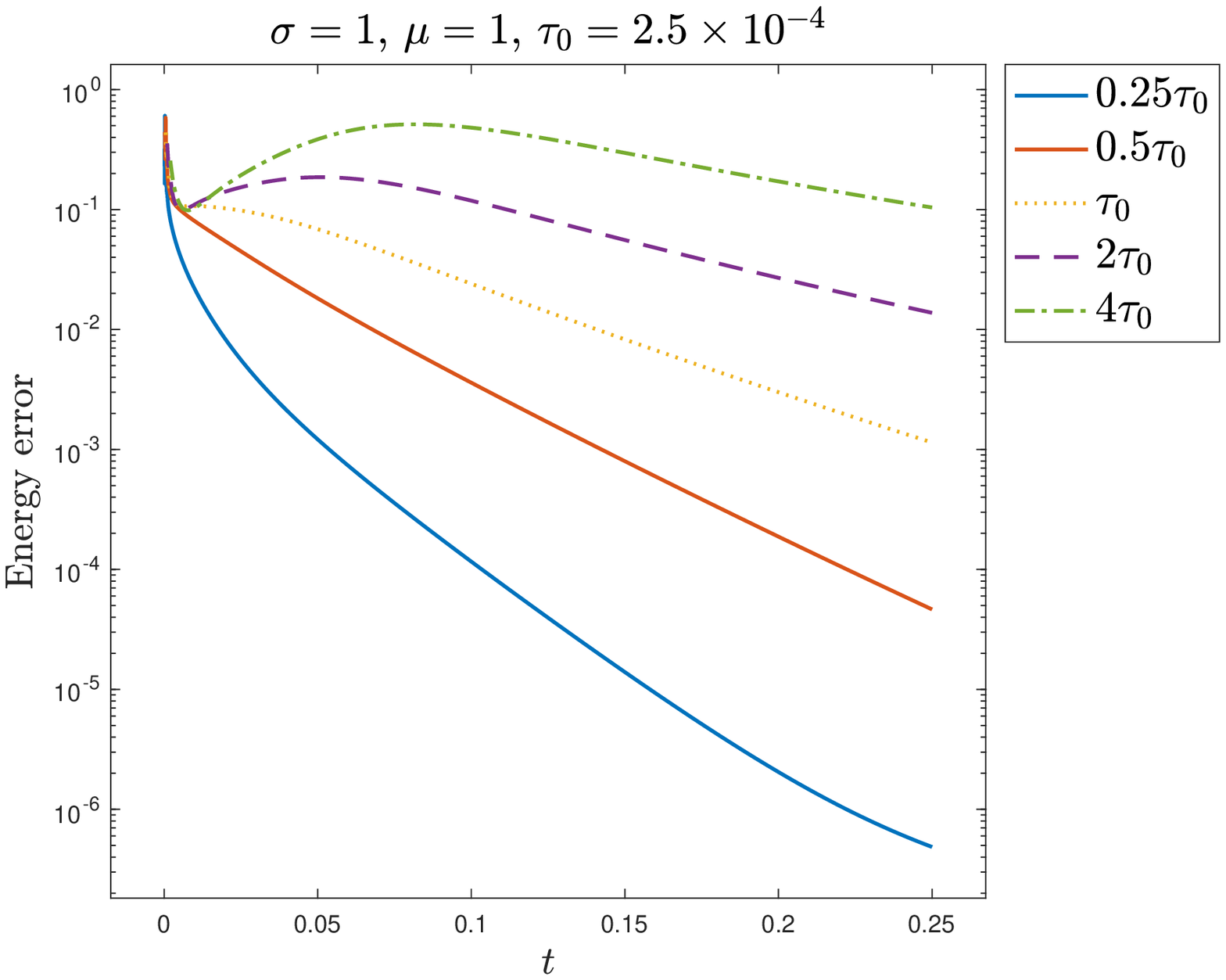}
\caption{History of energy error (using $1+9$ splitting). Left: $\sigma  \in \{ 1.0, 1.5 \}$ and $\mu \in \{ 1.0, 1.5\}$ with $\tau_0 = 2.5\times 10^{-4}$. Right: $\tau \in \{ 0.25 \tau_0, 0.5 \tau_0, \tau_0, 2\tau_0, 4\tau_0 \}$ with $(\mu, \sigma ) = (1,1)$ (Example \ref{exp3}).
}
\label{fig:his-exp3}
\end{figure}

\end{example}

\section{Conclusion} \label{sec:conclusion}
In this paper, we develop solution splitting algorithms for time-dependent multiscale problems. We construct a solution decomposition scheme based on the multiscale basis functions that represent the fine-scale details on a coarse grid. This solution decomposition technique can be combined with temporal splitting, and this leads to a smaller-dimensional matrix system to be solved in the temporal evolution. Moreover, this allows splitting fast and slow dynamics based on spatial resolution represented with multiscale basis functions. One can show that the whole splitting scheme for the spatial-discretized ODE system is consistent and unconditionally stable under appropriate conditions for the user-defined parameters. Numerical results are provided in order to illustrate the efficiency of the proposed method. We compare the results of splitting algorithms with backward Euler scheme. In the future, we will explore the application of such a decomposition method for more complicated cases such as unsteady parameter-dependent evolutionary problem. 

\section*{Acknowledgement}
YE would like to thank the partial support from NSF 1620318. YE would also like to acknowledge the support of Mega-grant of the Russian Federation Government (N 14.Y26.31.0013)

\bibliographystyle{abbrv}
\bibliography{references,references4,references1,references2,references3,decSol}
\end{document}